\newcommand{\siamyesno}[2]{#2}   
\siamyesno{
\documentclass[onetabnum,onefignum,nohypdvips,final]{siamart171218}
\usepackage{amssymb,amsmath,epsfig,verbatim,enumitem}
\newtheorem{remark}[theorem]{Remark}
\newtheorem{ass}[theorem]{Assumption}
}{
\documentclass[11pt,a4paper]{article}
\usepackage{amssymb,amsmath,amsthm,epsfig,verbatim,xcolor,enumitem}
\newtheorem{theorem}{\sc Theorem.}[section]
\newtheorem{lemma}[theorem]{\sc Lemma.}
\newtheorem{remark}[theorem]{\sc Remark.}

\newenvironment{AMS}%
{{\upshape\bfseries AMS subject classifications. }\ignorespaces}{}
\newenvironment{keywords}{{\upshape\bfseries Key words. }\ignorespaces}{}

}
\usepackage[normalem]{ulem}
\usepackage{ifpdf}  
\ifpdf
  \DeclareGraphicsExtensions{.eps,.pdf,.png,.jpg}
\else
  \DeclareGraphicsExtensions{.eps,.ps}
\fi
\newcommand{\bRplus}{{\bR}_{>0}}

\newcommand{\RZ}{{\bR} \slash {\mathbb Z}}
\newcommand{\bR}{{\mathbb R}}

\newcommand{\bN}{{\mathbb N}}

\newcommand{\drho}{\;{\rm d}\rho}
\newcommand{\ds}{\;{\rm d}s}

\newcommand{\Id}{I\!d}
\newcommand{\nabs}{\nabla_{\!s}}

\newcommand{\ratio}{{\mathfrak r}}
\newcommand{\Vh}{\underline{V}^h}

\def\epsilon{\varepsilon} 
\def\hat{\widehat}

\makeatletter
\renewcommand{\uuline}{%
  \bgroup
  \UL@setULdepth
  \markoverwith{%
    \lower\ULdepth\hbox{%
      \kern-.03em%
      \vtop{%
        \hrule width.2em%
        \kern 0.6pt 
        \hrule
      }%
      \kern-.03em%
    }%
  }%
  \ULon
}
\makeatother
\siamyesno{}{
\textwidth 455pt \oddsidemargin 0pt \evensidemargin 0pt \headsep
0pt \headheight 0pt \textheight 655pt \parskip 1pt \parindent 0pt

}

\begin{document}
\title{
Second order in time finite element schemes for \\
curve shortening flow and curve diffusion}
\author{Klaus Deckelnick\footnotemark[2]\ \and 
        Robert N\"urnberg\footnotemark[3]}

\renewcommand{\thefootnote}{\fnsymbol{footnote}}
\footnotetext[2]{Institut f\"ur Analysis und Numerik,
Otto-von-Guericke-Universit\"at Magdeburg, 39106 Magdeburg, Germany \\
{\tt klaus.deckelnick@ovgu.de}}
\footnotetext[3]{Dipartimento di Mathematica, Universit\`a di Trento,
38123 Trento, Italy \\ {\tt robert.nurnberg@unitn.it}}

\date{}

\maketitle

\begin{abstract}
We prove optimal error bounds for a second order in time finite element approximation of curve shortening flow in possibly higher codimension. In addition, we introduce a second order in time method for curve diffusion. Both schemes are based on variational formulations of strictly parabolic systems of partial differential equations that feature a tangential velocity which under discretization is beneficial for the mesh quality. In each time step only two linear systems need to be solved. Numerical experiments demonstrate second order convergence as well as asymptotic equidistribution.
\end{abstract} 

\begin{keywords} 
curve shortening flow; curve diffusion; finite elements; error analysis;
tangential motion; arbitrary codimension
\end{keywords}

\begin{AMS} 
65M60, 
65M12, 
65M15, 
35K55  
\end{AMS}

\renewcommand{\thefootnote}{\arabic{footnote}}

\siamyesno{
\pagestyle{myheadings}
\thispagestyle{plain}
\markboth{K. DECKELNICK AND R. N\"URNBERG}
{SECOND ORDER IN TIME FINITE ELEMENT SCHEMES}
}{}

\setcounter{equation}{0}
\section{Introduction} 

The aim of this paper is to introduce and analyse second order in time finite element schemes for the numerical solution of gradient
flows of the length functional
\begin{displaymath}
L(\Gamma) = \int_\Gamma 1 \ds.
\end{displaymath}
Here, $\Gamma$ denotes a closed curve in $\bR^d$ ($d \geq 2$) and ${\rm d}s$ is the arclength element. We consider both the $L^2$--gradient flow of $L$, i.e.\
the well--known curve shortening flow
\begin{equation} \label{eq:csf}
\mathcal V_{nor} = \varkappa \quad \mbox{ on } \Gamma(t),
\end{equation}
and the $H^{-1}$--gradient flow, which leads to curve diffusion
\begin{equation} \label{eq:cd}
\mathcal V_{nor} = - \nabs^2 \varkappa \quad \mbox{ on } \Gamma(t).
\end{equation}
We refer to \cite{GageH86,Grayson87} and 
\cite{TaylorC94,ElliottG97a,EscherMS98,GigaI98,GigaI99,DziukKS02,curves3d}
for some background on these two flows.
In the above, $(\Gamma(t))_{t \in [0,T]}$ denotes the family of evolving curves,  $\mathcal V_{nor}$ is the vector of normal velocities and $\varkappa$ the curvature vector 
of $\Gamma(t)$. Furthermore, $\nabs \phi$ is the vector of normal components of $\phi_s = \partial_s \phi$ for a vectorfield $\phi$.

In what follows we shall focus on a parametric description of the evolving curves, i.e.\ $\Gamma(t)= x(I,t)$ for some mapping $x:I \times [0,T] \to \bR^d$, where
$I = \RZ$ is the periodic unit interval. Note that \eqref{eq:csf} and \eqref{eq:cd} only prescribe the normal velocities, and that choosing a suitable tangential velocity component gives
rise to a (possibly degenerate) parabolic system of second and fourth order, respectively, for the position vector $x$.  Finite element schemes that discretize systems that are based on
reparameterizations obtained via the DeTurck trick have been proposed and analysed in \cite{DeckelnickD95} and \cite{ElliottF17} for the curve shortening flow and in \cite{cd}
for curve diffusion. In a series of papers, Barrett, Garcke and N\"urnberg developed schemes (the so--called BGN schemes) that equidistribute mesh points at the discrete level and used their
approach in particular to solve second and fourth order geometric evolution equations for curves, see \cite{triplej,triplejMC,curves3d,fdfi,bgnreview}. A different approach to enforce a uniform distribution
of mesh points along the curve has been proposed by Mikula and {\v{S}}ev{\v{c}}ovi{\v{c}} in \cite{MikulaS99}.  In all of the above schemes a backward Euler method is used for 
time discretization resulting in schemes of first order accuracy. In order to improve the efficiency and accuracy of the methods it is
desirable to develop schemes that are higher order in time and at the same time  enjoy good mesh and stability properties. 
In \cite{BalazovjechM11}, Balazovjech and Mikula suggest
a Crank--Nicolson type scheme for the curve shortening flow which turns out to be exact on shrinking circles. Mackenzie et al.\ propose in \cite{MackenzieNRI19}  a finite difference scheme
that uses a tangential mesh velocity in order to equidistribute grid points and which is second order accurate in time. In \cite{DuanLZ21} a high--order method is proposed for general
geometric functionals, which is based on a discretization of the first variation of the discrete energy and which uses a collocation scheme in time. All of the above schemes require the solution of a nonlinear
problem in each time step. More recently, higher order methods for both second and fourth order geometric flows based on the BGN--approach have been introduced by Jiang et al.\ in 
 \cite{JiangSZ24} using a Crank--Nicolson leapfrog scheme, and in \cite{JiangSZ24a} using backward differentiation formulae. The same authors propose in
 \cite{JiangSZZ25} a predictor--corrector scheme for surface diffusion, again on the basis of the BGN--discretization. The order of accuracy is verified numerically
 in terms of a suitably defined manifold distance. The above methods are applied in
 \cite{GarckeJSZ24preprint} to curve diffusion with a particular emphasis on ensuring that the schemes decrease length and conserve the enclosed area at the
 discrete level. 

 The aim of this paper is to introduce and analyse finite element schemes for both the curve shortening flow and curve diffusion in $\bR^d$, $d\geq2$,
that are based on systems obtained
 via the DeTurck trick. The corresponding systems have been introduced in \cite{DeckelnickD95} and \cite{cd}, respectively, and will be recalled in Section~\ref{sec:prelim}.
  The tangential velocity introduced in this way has the effect that the resulting systems are strictly parabolic and that the
 corresponding finite element discretizations are 
 well--behaved with respect to the distribution of grid points. In order to obtain schemes that are of second order in time, we shall adapt the predictor--corrector approach from \cite{JiangSZZ25}
 to this setting. More precisely, given the discrete solution at time $t_m$ and a time step size $\Delta t$, 
  the  predictor step computes an approximation at time $t_m + \frac12 \Delta t$ with the help of a known first order scheme, which is subsequently used in a Crank--Nicolson--type 
  step to calculate the discrete solution at time $t_m+\Delta t$. This idea will be carried out in Section~\ref{sec:csf} for the case of the curve shortening flow, extending the first order approach presented in \cite{DeckelnickD95}. 
We obtain a very simple method
  that only requires the solution of two linear systems in 
  each time step. Combining this new time discretization with a spatial discretization by continuous, piecewise linear finite elements, we shall prove that the resulting method is stable and convergent of optimal 
  order in space and time in the sense that
  \begin{itemize}
  \setlength\itemsep{0pt}
  \item the Dirichlet energy of the discrete solution is decreasing;
  \item the $L^2$--error is of order $O(h^2+(\Delta t)^2)$, the $H^1$--error is of order $O(h+(\Delta t)^2)$.
  \end{itemize}
Here we observe that our convergence proof relies on the time step size
being of order $O(h^\frac14)$, see Theorem~\ref{thm:main}, below, for details.
To the best of our knowledge, the rigorous proof of second order in time convergence in this setting is new in the literature. 
Let us point out that
recently second order in time schemes for mean curvature flow of axisymmetric
genus-1 surfaces have been proposed in \cite{LiWW25preprint}. 
In this case the profile curve evolves
according to an evolution law that coincides with the curve shortening
flow up to a lower order term. For an approach based on DeTurck's trick and
either a Crank-Nicolson scheme or a BDF2 method the authors obtain optimal
error bounds, but a discrete stability result does not seem to be available. 
Error estimates
for fully discrete schemes using backward differentiation formulae 
in the case of general parametric surfaces
have been obtained in \cite{KovacsLL19,BinzK23}. 
In these papers
the evolution is purely in normal direction, and again a proof of
discrete stability does not seem to be possible.
In addition, we refer to \cite{EsedogluG22} for a second-order accurate and
monotone threshold dynamics algorithm for a level set formulation of 
curve shorting flow.

Next, in Section~\ref{sec:cd} we shall apply the above ideas also  to curve diffusion, where we rely on the approach recently developed in \cite{cd}. Here, the discretization uses a splitting of the
fourth order problem into two second order problems by introducing a second variable in addition to the position vector. Special care needs to be taken with regard to the
relation between these two variables at the intermediate time step. We shall propose a finite element scheme that requires the solution of two linear systems in each step
and which is stable in the same sense as for the curve shortening flow. As the problem is much more complex, we are however not able to prove corresponding error bounds.
The numerical experiments which we present in Section~\ref{sec:nr} confirm
the error bounds in the case of the curve shortening flow, and suggest
that our new scheme for curve diffusion is also 
second order in time.

\setcounter{equation}{0}
\section{Preliminaries}  \label{sec:prelim}

\subsection{Notation} 
We denote the norm of the Sobolev space 
$W^{\ell,p}(I)$ ($\ell \in \bN_0$, $p \in [1, \infty]$)
by $\|\cdot \|_{\ell,p}$ and the 
semi-norm by $|\cdot |_{\ell,p}$. For
$p=2$, $W^{\ell,2}(I)$ will be denoted by
$H^{\ell}(I)$ with the associated norm and semi-norm written as,
respectively, $\|\cdot\|_\ell$ and $|\cdot|_\ell$.
The above are naturally extended to vector functions, and we will write 
$[W^{\ell,p}(I)]^d$ for a vector function with $d$ components.
In addition, throughout $c$ denotes a generic positive constant independent of 
the mesh parameter $h$ and the time step size $\Delta t$.
At times $\epsilon$ will play the role of a (small)
positive parameter, with $c_\epsilon>0$ depending on $\epsilon$, but
independent of $h$ and $\Delta t$.

\subsection{Curve shortening flow and curve diffusion}

Consider a family \linebreak $(\Gamma(t))_{t\in[0,T]}$ of evolving closed curves that are given by $\Gamma(t) = x(I,t)$, where
$ x : I \times [0,T] \to \bR^d$ satisfies $| x_\rho| > 0$ in $I \times [0,T]$. 
Then the unit tangent,
the curvature vector, the 
orthogonal projection onto the normal space and the 
vector of normal velocities of $\Gamma(t)$ are given by the 
following identities in $I$, see e.g.\ \cite{bgnreview}:
\begin{equation*} 
\tau = x_s = \frac{x_\rho}{| x_\rho |}, \quad
\varkappa =\frac{\tau_\rho}{| x_\rho |} =
\frac{1}{| x_\rho |}\left(\frac{x_\rho}{| x_\rho |} \right)_\rho, \quad 
P=\Id-\tau \otimes \tau, \quad
\mathcal V_{nor} = P x_t.
\end{equation*}
It is easily seen that if $x$ solves the system 
\begin{equation} \label{eq:DD95}
| x_\rho |^2 x_t = x_{\rho \rho} \quad \mbox{ in } I \times (0,T], \quad x(\cdot,0)= x_0,
\end{equation}
then $\Gamma(t)=x(I,t)$ is a solution of \eqref{eq:csf}. The system \eqref{eq:DD95} was used in \cite{DeckelnickD95} in order to derive a finite element scheme
for the approximation of curve shortening flow 
which allows for a simple error analysis. 

An extension of these ideas to curve diffusion was recently proposed in \cite{cd}. To do so, we introduce as a second variable $y=\frac{x_{\rho \rho}}{ | x_\rho |^2}$ and
consider the following system
\begin{equation} \label{eq:xt2}
|x_\rho|^2 x_t = - y_{\rho\rho} + F_{cd}(x_\rho,y,y_\rho) y
\quad \mbox{ in } I \times (0,T], \quad x(\cdot,0)= x_0,
\end{equation}
where $F_{cd}(a,b,c) \in \bR^{d \times d}$ is given by $F_{cd}=F_1+F_2$ with
\begin{subequations} 
\begin{align}
F_1(a,b,c) &= \bigl( 2 a \cdot c + | a |^2 | b |^2  \bigr)\Id, \label{eq:defF1} \\
F_2(a,b,c) &= 2 \bigl( c \otimes a -  a \otimes c \bigr) + 2 a \cdot b \bigl( a \otimes b - b \otimes a \bigr). \label{eq:defF2}
\end{align}
\end{subequations}
It is shown in \cite[\S 2.1]{cd} that if $x$ satisfies \eqref{eq:xt2}, then $\Gamma(t)=x(I,t)$ evolves according to \eqref{eq:cd} and satisfies in addition that 
$\frac{d}{dt} \int_I | x_\rho(\cdot,t) |^2 \leq 0$. 
A discrete analogue of this property also holds for the solution of a naturally defined finite element scheme, 
see \eqref{eq:cdfea} below, for which an optimal error bound
in the semidiscrete case is shown in \cite[Theorem~3.1]{cd}.

\subsection{Finite elements and finite differences} 

Let $[0,1]=\bigcup_{j=1}^J I_j$, $J\geq3$, be a
decomposition of $[0,1]$ into the intervals 
$I_j=[\rho_{j-1},\rho_j]$, where $0=\rho_0 < \rho_1 < \ldots < \rho_J=1$ and 
set $\mathcal G^h= \{\rho_1,\ldots, \rho_J \}$.  Due to the periodicity of
$\RZ$ we identify $\rho_0 = 0$ with $\rho_J=1$ and for convenience set
$\rho_{J+1} = \rho_1$ and $\rho_{-1} = \rho_{J-1}$.
Let us define $h_j= \rho_j - \rho_{j-1}$, the grid size $h=
\max_{1 \leq j \leq J} h_j$ and assume in what follows that
\begin{equation} \label{eq:decass}
|  h_j - h_{j-1} | \leq c h^2, \quad h \leq c h_j,  \qquad j=1,\ldots,J.
\end{equation}
{From} now on we let $(\cdot,\cdot)$ denote the $L^2$--inner product on $I$.
For two piecewise continuous functions, with possible jumps at the 
nodes $\{\rho_j\}_{j=1}^J$, we define the mass lumped $L^2$--inner product 
\begin{equation*} 
( u, v )^h = \tfrac12\sum_{j=1}^J h_j
\left[(u\cdot v)(\rho_j^-) + (u \cdot v)(\rho_{j-1}^+)\right],
\end{equation*}
where $(u \cdot v)(\rho_j^\pm)=\underset{\delta\searrow 0}{\lim}\ 
(u \cdot v)(\rho_j\pm\delta)$. 
We define the finite element spaces
\[
V^h = \{\chi \in C^0(I) : \chi\!\mid_{I_j} 
\text{ is affine},\ j=1,\ldots, J\} \quad \text{and}\quad \Vh = [V^h]^d,
\]
and let $I_h : C^0(I) \to V^h$ denote the standard interpolation
operator at the nodes $\{\rho_j\}_{j=1}^J$. We use the same notation for the
interpolation of vector-valued functions.
It is well-known that for 
$k \in \{ 0,1 \}$, $\ell \in \{ 1,2 \}$ and $p \in [2,\infty]$ it holds that
\begin{subequations}
\begin{alignat}{2}
h^{\frac 1p - \frac 1r} \| \eta_h \|_{0,r} 
+ h | \eta_h |_{1,p} & \leq C \| \eta_h \|_{0,p} 
\qquad && \forall\ \eta_h \in V^h, \qquad r \in [p,\infty], 
\label{eq:inverse} \\
| \eta - I_h \eta |_{k,p} & \leq Ch^{\ell-k} | \eta |_{\ell,p} 
\qquad && \forall\ \eta \in W^{\ell,p}(I). \label{eq:estIh} 
\end{alignat}
\end{subequations}
For later use we also note that for $j=1,\ldots,J$
\begin{equation} \label{eq:inteq}
\int_{I_j} (u - I_h u)_\rho\, \eta_{h,\rho} \drho = 0 \qquad \forall\ u \in H^1(I),
\quad \forall\ \eta_h \in V^h,
\end{equation}
and, as a consequence of the error formula for the trapezium rule,
\begin{align} \label{eq:trapezrule}
&\int_{I_j} \eta_h \chi_h \drho 
- \tfrac12 h_j \left[(\eta_h \chi_h)(\rho_j) + (\eta_h \chi_h)(\rho_{j-1})\right]
= - \tfrac16\,h_j^2 \int_{I_j} \, \eta_{h,\rho} \, \chi_{h,\rho} \drho 
\nonumber \\ & \hspace{9cm} \forall\ \eta_h, \chi_h \in V^h.
\end{align}
Even though we will formulate our numerical method as a finite element scheme, it will be
convenient to carry out part of the error analysis in a finite difference setting. With a 
grid function $v: \mathcal G_h \to \bR^d$ we associate the values
$v_j=v(\rho_j)$ and introduce the finite difference operators
\begin{displaymath}
\delta^- v_j= \frac{v_j - v_{j-1}}{h_j}, \quad 
\delta^+ v_j= \frac{v_{j+1} - v_j}{h_{j+1}}, \quad 
\delta^2 v_j= \frac{\delta^+v_j - \delta^- v_j}{\frac12(h_j+h_{j+1})}.
\end{displaymath}
For two grid functions $v,w: \mathcal G_h \to \bR^d$ one has the following summation by parts formula:
\begin{equation} \label{eq:sbp}
\sum_{j=1}^J h_j \, \delta^- v_j \cdot \delta^- w_j = - \sum_{j=1}^J \tfrac12 (h_j+h_{j+1}) \, v_j \cdot \delta^2 w_j.
\end{equation}

In order to discretize in time, let $t_m=m \Delta t$, $m=0,\ldots,M$, 
with the uniform time step $\Delta t = \frac TM >0$, and similarly
$t_{m+\frac12}= (m+\frac12) \Delta t$, $m=0,\ldots,M-1$, for later use.

Now the finite element approximation for \eqref{eq:DD95} from
\cite{DeckelnickD95} with mass lumping can be formulated as follows.
For $m \geq 0$, given $x^m_h \in \Vh$, find $x^{m+1}_h\in\Vh$ such that
\begin{equation} \label{eq:DD95:fea}
\left(
\frac{x^{m+1}_h-x^m_h}{\Delta t} , \eta_h |x^m_{h,\rho}|^2 \right)^h
+ \left( x^{m+1}_{h,\rho} , \eta_{h,\rho} \right) = 0
\qquad \forall\ \eta_h \in \Vh.
\end{equation}
Similarly, the scheme for \eqref{eq:xt2} from \cite{cd} is given by:
For $m \geq 0$, given $(x^m_h, y^m_h) \in \Vh \times \Vh$ find 
$(x^{m+1}_h, y^{m+1}_h) \in \Vh \times \Vh$ such that
\begin{subequations} \label{eq:cdfea}
\begin{align}
& \left( \frac{x^{m+1}_h-x^m_h}{\Delta t} , \chi_h |x^m_{h,\rho}|^2 \right) 
- \left( y^{m+1}_{h,\rho} , \chi_{h,\rho} \right) 
= 2 \left( (y^{m+1}_{h,\rho} \cdot x^m_{h,\rho}) y^m_h , \chi_h \right)
\nonumber \\ & \
+ \left( (y^m_h \cdot y^{m+1}_h) y^m_h , \chi_h | x^m_{h,\rho} |^2 \right)
+  \left( F_2(x^m_{h,\rho},y^m_h,y^m_{h,\rho})y^{m+1}_h , \chi_h \right) 
\quad \forall\ \chi_h \in \Vh, \label{eq:cdfeaa} \\
& \left( y^{m+1}_h , \eta_h |x^m_{h,\rho}|^2 \right) 
+ \left( x^{m+1}_{h,\rho} , \eta_{h,\rho} \right) = 0
\quad \forall\ \eta_h \in \Vh. \label{eq:cdfeab}
\end{align}
\end{subequations}

\setcounter{equation}{0}
\section{A second order scheme for the curve shortening flow} \label{sec:csf}

As already outlined in the introduction, we aim to apply the idea from \cite{JiangSZZ25} of using a predictor--corrector approach to the
first order method \eqref{eq:DD95:fea} from \cite{DeckelnickD95}. This leads to the following finite element scheme:
Let $x^0_h= I_h x_0 \in \Vh$. Then, for $m=0,\ldots,M-1$, 
first find $x^{m+\frac12}_h\in\Vh$ such that
\begin{subequations} \label{eq:fea}
\begin{equation} \label{eq:pred}
\left(
\frac{x^{m+\frac12}_h-x^m_h}{\frac12\Delta t} , \eta_h |x^m_{h,\rho}|^2 \right)^h
+ \left( x^{m+\frac12}_{h,\rho} , \eta_{h,\rho} \right) = 0
\qquad \forall\ \eta_h \in \Vh,
\end{equation}
and then find $x^{m+1}_h \in \Vh$ such that
\begin{equation} \label{eq:csfd}
\left( \frac{x^{m+1}_h-x^m_h}{\Delta t} , \eta_h 
|x^{m+\frac12}_{h,\rho}|^2 \right)^h
+ \tfrac12 \left(x^{m+1}_{h,\rho} + x^{m}_{h,\rho} , \eta_{h,\rho} \right) =  0
\qquad \forall\ \eta_h \in \Vh.
\end{equation}
\end{subequations}

With a mild assumption on the discrete length elements, the scheme
\eqref{eq:fea} has a unique solution that is unconditionally stable.

\begin{theorem} \label{thm:stab}
Assume that $|x^m_{h,\rho}| > 0$ in $I$.
Then there exist unique solutions 
$x^{m+\frac12}_h \in \Vh$ to \eqref{eq:pred}
and $x^{m+1}_h\in \Vh$ to \eqref{eq:csfd}.
Moreover, any solution to \eqref{eq:csfd} satisfies the stability estimate
\begin{equation} \label{eq:fdstab}
\tfrac12 |x^{m+1}_{h}|^2_1
+ \frac1{\Delta t} \left( | x^{m+1}_{h} - x^m_h|^2, 
|x^{m+\frac12}_{h,\rho} |^2 \right)^h
= \tfrac12 |x^m_{h}|^2_1 .
\end{equation}
\end{theorem}
\begin{proof}
The well-posedness of the linear system \eqref{eq:pred} follows from the unique
solvability of the homogeneous system. Let $X_h \in \Vh$ satisfy
\begin{equation} \label{eq:homo}
2\left(
\frac{X_h}{\Delta t} , \eta_h |x^m_{h,\rho}|^2 \right)^h
+ \left( X_{h,\rho} , \eta_{h,\rho} \right) = 0
\qquad \forall\ \eta_h \in \Vh.
\end{equation}
Choosing $\eta_h = X_h$ in \eqref{eq:homo} yields
\[
2 \left( |X_h|^2, |x^m_{h,\rho}|^2 \right)^h + \Delta t |X_h|_1^2 = 0.
\]
Hence $X_h = X \in \bR^d$ must be a constant satisfying 
$2|X|^2 (1,|x^m_{h,\rho}|^2 ) =0$. Therefore our assumptions imply that
$X=0$ and so \eqref{eq:pred} has a unique solution. 
In order to show the well-posedness of the linear system \eqref{eq:csfd} we 
proceed as before. 
In particular, we obtain that the solution $X_h$ of the corresponding 
homogeneous system is a constant $X \in \bR^d$ with
$2|X|^2 (1,|x^{m+\frac 12}_{h,\rho}|^2)=0$. If $X \neq 0$, then
$x^{m+\frac12}_{h,\rho}= 0$ in $I$ and hence 
$x^{m+\frac12}_{h} = x_c \in \bR^d$ is a constant. Then choosing 
$\eta_h = x_c - x^m_h$ in \eqref{eq:pred} implies $x^m_h = x_c$, which 
would contradict our assumption on $x^m_h$. Thus, $X=0$.

In order to prove \eqref{eq:fdstab}, we choose 
$\eta_h = x^{m+1}_h - x^m_h$ in \eqref{eq:csfd} to yield that 
\[
\frac1{\Delta t} \left( | x^{m+1}_{h} - x^m_h|^2, 
|x^{m+\frac12}_{h,\rho} |^2 \right)^h
+ \tfrac12 |x^{m+1}_{h}|^2_1
- \tfrac12 |x^m_{h}|^2_1 = 0,
\]
which is our desired result \eqref{eq:fdstab}. 
\end{proof}

\vspace{3mm}
Our main result is the following optimal error estimate.
\begin{theorem} \label{thm:main}
Suppose that \eqref{eq:DD95}
has a smooth solution on the time interval $[0,T+\delta)$ satisfying
\begin{equation} \label{eq:regul}
c_0 \leq | x_\rho | \leq C_0 \quad \mbox{ in } I \times [0,T]
\end{equation}
for some constants $\delta, c_0, C_0 \in \bRplus$.
Then there exist $\gamma > 0$ and $h_0 > 0$ such that if $0 < h \leq h_0$ and
$\Delta t \leq \gamma h^{\frac14}$,
then \eqref{eq:fea} has a unique 
solution $(x^{m}_h)_{m=0,\ldots,M}$, and the following error bounds hold:
\begin{subequations} 
\begin{align}
\max_{0 \leq m \leq M} \| x(\cdot,t_m) - x^m_h \|_0 
+ h \max_{0 \leq m \leq M} | x(\cdot,t_m) - x^m_h |_1 & \leq c \bigl( h^2 + (\Delta t)^2 \bigr), \label{eq:l2h1} \\
\sum_{m=0}^{M-1} \Delta t \, \left\| x_t(\cdot,t_{m+\frac12}) - \frac{x^{m+1}_h - x^m_h}{\Delta t} \right\|^2_0 & \leq c \bigl( h^4 + (\Delta t)^4 \bigr). \label{eq:timedif}
\end{align}
\end{subequations}
\end{theorem}
\begin{proof} 
In what follows we abbreviate $x^k = x(\cdot,t_k)$ for 
$k \in \{ m,m+\tfrac12 \}$, as well as
\begin{equation} \label{eq:eh}
e^k_h= x^k_h - I_h x^k \in \Vh, \quad d^k= x^k - I_h x^k, 
\qquad k \in \{ m,m+\tfrac12 \}.
\end{equation}
Furthermore, for $m \in \lbrace 0,1,\ldots,M \rbrace$ and 
$\gamma h_0^\frac14 \leq \delta$ we introduce the following error measure
\begin{align} \label{eq:defEm}
E^m &= \tfrac12 | e^m_{h} |^2_1   +2 \left(x_t^{m+\frac12} \cdot e^{m}_{h,\rho},x^{m+\frac12}_\rho \cdot d^{m+\frac12} \right) \nonumber \\ & \quad
+ \tfrac16 \sum_{j=1}^J h_j^2 \int_{I_j} \frac{x^{m+1}_\rho - x^m_\rho}{\Delta t}  \cdot e^m_{h,\rho} \, | x^{m+\frac12}_\rho |^2 \drho.
\end{align}
It follows from the interpolation estimate \eqref{eq:estIh} 
and the smoothness of $x$ that
\begin{equation} \label{eq:Eequiv}
\tfrac14 | e^m_{h} |^2_1 - c h^4 \leq E^m \leq \tfrac34 | e^m_{h} |^2_1 + c h^4,
\quad m = 0,\ldots,M.
\end{equation}
The specific form of $E^m$ will allow us to obtain a superconvergence result for $| e^m_h |_1$. More precisely, we shall
show by induction that $x^m_h\in\Vh$ exists uniquely and satisfies
\begin{equation} \label{eq:induction}
E^m \leq \bigl( h^4 + (\Delta t)^4 \bigr) e^{\mu t_m},\quad m = 0,\ldots,M,
\end{equation}
provided that $0 <h \leq h_0$, 
$\Delta t \leq \gamma h^{\frac14}$
and $\gamma>0, \mu>0$ are chosen appropriately. 
Clearly the assertion holds for $m=0$ since 
$x_h^0 = I_0 x_0$ and so $E^0=0$.
Suppose for some $0\leq m \leq M-1$ that $x^m_h\in\Vh$ exists uniquely and \eqref{eq:induction} is satisfied. 
Then we deduce from \eqref{eq:eh}, 
\eqref{eq:regul}, \eqref{eq:estIh}, \eqref{eq:inverse}, 
\eqref{eq:Eequiv}, \eqref{eq:induction} and the smoothness of $x$ that
\begin{align*}
| x^m_{h,\rho} | & \geq | x^m_\rho | - | d^m_\rho | - | e^m_{h,\rho} |
\geq c_0 - c h - c h^{-\frac12} | e^m_h |_1 \\ &
\geq c_0 - ch - c h^{-1/2} \bigl( 2 (E^m)^\frac12 + c h^2 \bigr) 
\geq c_0 - ch - c   h^{-\frac12} e^{\frac12 \mu T} \bigl( h^2+ (\Delta t)^2 \bigr) \\ & 
\geq c_0 - c h_0 e^{\frac12 \mu T} - c \gamma^2 e^{\frac12 \mu T}
\geq \tfrac12 c_0,
\end{align*}
provided that  $c h_0 e^{\frac12 \mu T} \leq \frac{c_0}{4}$ and $c \gamma^2 e^{\frac12 \mu T} \leq \frac{c_0}{4}$. Arguing in the same way in order to bound $| x^m_{h,\rho} |$ from above, and using \eqref{eq:estIh} for similar bounds
on $|(I_h x^m)_\rho|$, we obtain that
\begin{equation} \label{eq:xhbound1}
\tfrac12 c_0 \leq | (I_h x^m)_\rho | \leq 2 C_0 \quad\text{and}\quad
\tfrac12 c_0 \leq | x^m_{h,\rho} | \leq 2 C_0 \quad \mbox{ in } I.
\end{equation}
We shall choose $\mu$ at the end of the proof and this choice then
determines the smallness conditions for $h_0$ and $\gamma$. \\
Theorem~\ref{thm:stab} together with \eqref{eq:xhbound1} implies that $x^{m+\frac12}_h\in\Vh$ exists uniquely.
The following lemma provides an estimate for the discrete error
$e^{m+\frac12}_{h}$ arising from \eqref{eq:pred}. We use a reformulation of
\eqref{eq:pred} into a finite difference scheme to prove it.
\begin{lemma} \label{lem:predest} 
There exists $c>0$ such that 
\begin{equation} \label{eq:predest}
| e^{m+\frac12}_{h} |^2_1 
\leq c | e^m_h|^2_1 + c \bigl( h^4 +  (\Delta t)^4 \bigr).
\end{equation}
\end{lemma}
\begin{proof}
Let us begin by writing \eqref{eq:pred} as a finite difference scheme.
On defining
\begin{displaymath}
q^k_{h,j}= | \delta^- x^k_{h,j}| 
= \left| \frac{x^k_{h,j}-x^k_{h,j-1}}{h_j}\right|, \quad j=1,\ldots,J+1,
\end{displaymath}
for $k \in \{m, m+\frac12\}$, we see that \eqref{eq:pred} can be written as
\begin{displaymath}
\tfrac12\bigl( h_j (q^m_{h,j})^2 + h_{j+1} (q^m_{h,j+1})^2 \bigr) \frac{x^{m+\frac12}_{h,j} - x^m_{h,j}}{\frac12 \Delta t} - \bigl( \delta^+ x^{m+\frac12}_{h,j} - \delta^- x^{m+\frac12}_{h,j} \bigr)=0, \quad
j=1,\ldots,J,
\end{displaymath}
or equivalently
\begin{equation} \label{eq:fd}
x^{m+\frac12}_{h,j} - x^m_{h,j} - \Delta t \, \alpha^m_j \,   \delta^2 x^{m+\frac12}_{h,j}  = 0, \qquad j=1,\ldots,J,
\end{equation}
where
\begin{equation*} 
\alpha^m_j= \frac{\frac12 (h_j+h_{j+1})}{ h_j (q^m_{h,j})^2 + h_{j+1} (q^m_{h,j+1})^2}.
\end{equation*}

In order to write down the error relation we define the grid functions $e^k$
via
\begin{equation} \label{eq:ej}
e^k_j= x^k_{h,j} - x^k_j= x^k_{h,j} - (I_hx^k)_j=e^k_h(\rho_j),
\quad j=1,\ldots,J,
\end{equation}
for $k \in \{m, m+\frac12\}$, 
and then combine \eqref{eq:fd} with \eqref{eq:DD95} to obtain that
\begin{align*}
& e^{m+\frac12}_j - e^m_j -  \Delta t \, \alpha^m_j \,   \delta^2 e^{m+\frac12}_{j} \\ & \quad
= \tfrac12 \Delta t \left(  x_t(\rho_j,t_m) - \frac{x^{m+\frac12}_j-x^m_j}{\frac12 \Delta t} \right)+ \frac{\Delta t}{2 | x^m_\rho(\rho_j) |^2} \bigl( x^{m+\frac12}_{\rho \rho}(\rho_j) -x^m_{\rho \rho}(\rho_j) \bigr) \\ & \qquad 
+ \Delta t \left(  \alpha^m_j - \frac{1}{2 | x^m_\rho(\rho_j) |^2} \right) \delta^2 x^{m+\frac12}_j  + \Delta t \, \frac{1}{2 | x^m_\rho(\rho_j) |^2} \bigl( 
\delta^2 x^{m+\frac12}_j - x^{m+\frac12}_{\rho \rho}(\rho_j) \bigr) \\
& \quad
=:\sum_{\ell=1}^4 f^{m}_{\ell,j}, \quad
 j=1,\ldots,J.
\end{align*}
Let us multiply the above equation by $-\frac{h_j+h_{j+1}}{2} \, \delta^2 e^{m+\frac12}_j$, sum over $j=1,\ldots,J$ and use \eqref{eq:sbp}. Then
\begin{align}
& \sum_{j=1}^J h_j \,\delta^- (e^{m+\frac12}_j -e^m_j) \cdot \delta^- e^{m+\frac12}_j 
+  \Delta t  \sum_{j=1}^J \tfrac12 (h_j+h_{j+1}) \alpha^m_j \,   | \delta^2 e^{m+\frac12}_j |^2 \nonumber \\ & \qquad\qquad
= - \sum_{\ell=1}^4 \sum_{j=1}^J \tfrac12(h_j+h_{j+1}) \,   f^{m}_{\ell,j} \cdot \delta^2 e^{m+\frac12}_j. \label{eq:erra0}
\end{align}
Note that \eqref{eq:xhbound1} implies that 
\begin{equation} \label{eq:qbound}
\tfrac12 c_0 \leq q^m_{h,j} \leq 2 C_0,\qquad j = 1,\ldots,J,
\end{equation}
and hence in particular
\begin{displaymath}
\alpha^m_j \geq \tfrac12 \frac{1}{(q^m_{h,j})^2 + (q^m_{h,j+1})^2} \geq \frac{1}{16 C_0^2}.
\end{displaymath}
Using this bound in \eqref{eq:erra0}, together with the relation $(a-b) \cdot a \geq \frac{1}{2} | a|^2 - \frac{1}{2} |b|^2$, we derive
\begin{align} 
& \tfrac12  \sum_{j=1}^J h_j  | \delta^- e^{m+\frac12}_j |^2 
- \tfrac12  \sum_{j=1}^J h_j  | \delta^- e^m_j |^2 + \frac{\Delta t}{16 C_0^2} \sum_{j=1}^J \tfrac12(h_j+h_{j+1})  | \delta^2 e^{m+\frac12}_j |^2 \nonumber \\
& \qquad\qquad
\leq - \sum_{\ell=1}^4 \sum_{j=1}^J \tfrac12(h_j+h_{j+1}) f^{m}_{\ell,j} \cdot \delta^2 e^{m+\frac12}_j. \label{eq:erra1}
\end{align}
In order to deal with the first term on the right hand side of \eqref{eq:erra1} we use \eqref{eq:sbp} and obtain
\begin{displaymath}
- \sum_{j=1}^J \tfrac12(h_j+h_{j+1}) f^{m}_{1,j} \cdot \delta^2 e^{m+\frac12}_j =  \sum_{j=1}^J h_j \, \delta^- f^{m}_{1,j} \cdot \delta^- e^{m+\frac12}_j.
\end{displaymath}
A short calculation shows that
\begin{align*}
| \delta^- f^{m}_{1,j} | & = \left| \int_{t_m}^{t_{m+\frac12}} \frac{1}{h_j} \int_{\rho_{j-1}}^{\rho_j} x_{t \rho}(\sigma,t_m) - x_{t \rho}(\sigma,\tau)\; {\rm d}\sigma {\rm d}\tau \right| \\ & 
\leq c (\Delta t)^2 \max_{t\in[0,T+\frac{\delta}{2}]}\Vert x_{tt,\rho}(\cdot,t) \Vert_{0,\infty},
\end{align*}
so that
\begin{align} \label{eq:f1}
-  \sum_{j=1}^J \tfrac12(h_j+h_{j+1}) f^{m}_{1,j} \cdot \delta^2 e^{m+\frac12}_j & \leq c (\Delta t)^2 \left( \sum_{j=1}^J h_j | \delta^- e^{m+\frac12}_j |^2 \right)^{\frac12} \nonumber \\ &
\leq \tfrac18 \sum_{j=1}^J h_j | \delta^- e^{m+\frac12}_j |^2 + c
 (\Delta t)^4.
\end{align}
Observing that
\begin{displaymath} 
| \delta^- f^m_{2,j} | \leq c (\Delta t)^2 \left( \max_{t\in [0,T+\frac{\delta}{2}]} \Vert x_{t,\rho \rho \rho}(\cdot,t) \Vert_{0,\infty} 
+ \max_{t\in [0,T+\frac{\delta}{2}]} \Vert x_{t,\rho \rho}(\cdot,t) \Vert_{0,\infty} \right),
\end{displaymath}
the second term can be treated with the same argument, so that
\begin{equation} \label{eq:f2}
-  \sum_{j=1}^J \tfrac12(h_j+h_{j+1}) f^{m}_{2,j} \cdot \delta^2 e^{m+\frac12}_j \leq \tfrac18  \sum_{j=1}^J h_j | \delta^- e^{m+\frac12}_j |^2 + c (\Delta t)^4.
\end{equation}
Next, let $q^m_j = |\delta^- x^m_j| = |(I_h x^m)_{\rho|_{I_j}}|$, for which we
have $\frac12c_0 \leq q^m_j \leq 2 C_0$ in view of \eqref{eq:xhbound1}. 
Furthermore, Taylor expansions together with \eqref{eq:decass} imply
\begin{subequations} 
\begin{align}
\left|  (q^m_j)^2 + (q^m_{j+1})^2 - 2 | x^m_\rho(\rho_j) |^2  \right| &\leq c h^2, \qquad j=1,\ldots,J, \label{eq:taylor1} \\
|\delta^2 x^{m+\frac12}_j - x^{m+\frac12}_{\rho \rho}(\rho_j)| & \leq c h^2,
\qquad j=1,\ldots,J. \label{eq:taylor2}
\end{align}
\end{subequations}
Then we use \eqref{eq:decass},  \eqref{eq:taylor1} and \eqref{eq:qbound} to estimate
\begin{align*}
&\left|  \alpha^m_j -  \frac{1}{ 2 | x^m_\rho(\rho_j) |^2} \right| \leq \left|  \frac{\frac12 (h_j+h_{j+1})}{ h_j (q^m_{h,j})^2 + h_{j+1} (q^m_{h,j+1})^2} - \frac{1}{ (q^m_{h,j})^2 + (q^m_{h,j+1})^2}  \right| \\
& \quad + \left| \frac{1}{ (q^m_{h,j})^2 + (q^m_{h,j+1})^2} -  \frac{1}{ (q^m_{j})^2 + (q^m_{j+1})^2} \right| +\left|   \frac{1}{ (q^m_{j})^2 + (q^m_{j+1})^2} - \frac{1}{ 2 | x^m_\rho(\rho_j) |^2} \right| \\
& \quad \leq  c h^{-1} | h_{j+1} - h_j | \, | (q^m_{h,j})^2 - (q^m_{h,j+1})^2 | 
+ c \bigl( | q^m_{h,j} - q^m_j | + | q^m_{h,j+1} - q^m_{j+1} | \bigr) \\ & \quad \quad  + c | (q^m_j)^2 + (q^m_{j+1})^2 - 2 | x^m_\rho(\rho_j) |^2 |  \\ 
& \quad \leq ch | q^m_{j+1} - q^m_j | + c \bigl( | \delta^- e^m_j | + | \delta^- e^m_{j+1} | \bigr) + c h^2 \\
& \quad
\leq c \bigl( | \delta^- e^m_j | + | \delta^- e^m_{j+1} | \bigr) + c h^2.
\end{align*}
Combining this bound with \eqref{eq:qbound} and \eqref{eq:decass}, we obtain 
\begin{align}
&
- \sum_{j=1}^J \tfrac12(h_j+h_{j+1}) f^{m}_{3,j} \cdot \delta^2 e^{m+\frac12}_j 
\nonumber \\ & \quad
 \leq c \Delta t \sum_{j=1}^J \tfrac12(h_j+h_{j+1}) \bigl(  | \delta^- e^m_j | + | \delta^- e^m_{j+1} | + h^2 \bigr) | \delta^2 e^{m+\frac12}_j |
\nonumber \\ & \quad
\leq c \Delta t \left( \left(  \sum_{j=1}^J h_j  | \delta^- e^m_j |^2 \right)^{\frac12} + \left(  \sum_{j=1}^J h_{j+1}  | \delta^- e^m_{j+1} |^2 \right)^{\frac12} + h^2 \right) \cdot \nonumber \\ & \qquad \qquad \cdot
\left(  \sum_{j=1}^J \tfrac12(h_j+h_{j+1}) | \delta^2 e^{m+\frac12}_j |^2 \right)^{\frac12}  \nonumber \\ & \quad
\leq \frac{\Delta t}{32 C_0^2}  \sum_{j=1}^J \tfrac12(h_j+h_{j+1}) | \delta^2 e^{m+\frac12}_j |^2 + c \Delta t \left(  \sum_{j=1}^J h_j | \delta^- e^m_j |^2 + h^4 \right). \label{eq:f3}
\end{align}
Finally, a similar argument together with \eqref{eq:taylor2} yields
\begin{equation} \label{eq:f4}
- \sum_{j=1}^J \tfrac12(h_j+h_{j+1}) f^{m}_{4,j} \cdot \delta^2 e^{m+\frac12}_j \leq \frac{\Delta t}{32 C_0^2}  \sum_{j=1}^J \tfrac12(h_j+h_{j+1}) | \delta^2 e^{m+\frac12}_j |^2 + c h^4 \Delta t.
\end{equation}
If we insert \eqref{eq:f1}, \eqref{eq:f2}, \eqref{eq:f3} and \eqref{eq:f4} into \eqref{eq:erra1} we derive
\begin{displaymath}
 \sum_{j=1}^J h_j   |  \delta^- e^{m+\frac12}_j |^2 \leq c  \sum_{j=1}^J h_j  | \delta^- e^{m}_j |^2+ c \bigl( h^4 + (\Delta t)^4 \bigr).
\end{displaymath}
Combining this with the fact that
\[
 \sum_{j=1}^J h_j | \delta^-e^k_j|^2 = 
  \| e^k_{h,\rho} \|^2_0
= | e^k_{h} |^2_1,
\]
recall \eqref{eq:ej} and \eqref{eq:eh}, yields the desired result
\eqref{eq:predest}. 
\end{proof}

\vspace{3mm}
With the help of Lemma~\ref{lem:predest} and \eqref{eq:Eequiv} we can argue in a very similar way 
as for \eqref{eq:xhbound1} to show that
\begin{equation} \label{eq:xhbound2}
\tfrac12 c_0 \leq | (I_h x^{m+\frac12})_\rho | \leq 2 C_0 \quad\text{and}\quad
\tfrac12 c_0 \leq | x^{m+\frac12}_{h,\rho} | \leq 2 C_0 \quad \mbox{ in } I,
\end{equation}
provided that $c h_0 e^{\frac12 \mu T} \leq \frac{c_0}{4}$ and $c \gamma^2 e^{\frac12 \mu T} \leq \frac{c_0}{4}$ with a possibly larger constant $c$ compared to
\eqref{eq:xhbound1}. 
Hence we may infer from Theorem~\ref{thm:stab} that the solution
$x^{m+1}\in\Vh$ of \eqref{eq:csfd} exists uniquely. As a second step in our error analysis we prove

\begin{lemma} 
There exists $c>0$ such that
\begin{equation} \label{eq:lem2}
\frac{c_0^2}{8 \Delta t} \| e^{m+1}_h - e^m_h \|^2_0 +  E^{m+1} - E^m 
 \leq c \Delta t | e^{m+\frac12}_{h} |^2_1 + c \Delta t | e^{m+1}_{h} |^2_1 + c \Delta t ( h^4 + (\Delta t)^4 ).
\end{equation} 
\end{lemma}
\begin{proof} 
We deduce from
 \eqref{eq:DD95} that
\begin{displaymath}
\left( x_t^{m+\frac12} , \eta_h\, | x_\rho^{m+\frac12} |^2 \right)
+ \left(x^{m+\frac12}_\rho , \eta_{h,\rho} \right) =0 \qquad \forall\ \eta_h\in
\Vh,
\end{displaymath}
while \eqref{eq:csfd} along with \eqref{eq:trapezrule} implies
\begin{align*}
& \left( \frac{x^{m+1}_h-x^m_h}{\Delta t} , \eta_h 
|x^{m+\frac12}_{h,\rho}|^2 \right) \\ & \quad
+ \tfrac16 \sum_{j=1}^J h_j^2 \int_{I_j} \frac{x^{m+1}_{h,\rho}-x^m_{h,\rho}}{\Delta t} \cdot \eta_{h,\rho} \, |x^{m+\frac12}_{h,\rho}|^2 \drho
+ \tfrac12 \left(x^{m+1}_{h,\rho} + x^{m}_{h,\rho} , \eta_{h,\rho} \right) =  0
\end{align*}
for all $\eta_h \in \Vh$. Combining these equations with \eqref{eq:inteq} 
yields the following error relation
\begin{align}
&
\left( \frac{e^{m+1}_h- e^m_h}{\Delta t} , \eta_h\, | x^{m+\frac12}_{h,\rho} |^2
\right) + \tfrac16 \sum_{j=1}^J h_j^2 \int_{I_j} \frac{e^{m+1}_{h,\rho}-e^m_{h,\rho}}{\Delta t} \cdot \eta_{h,\rho} \, |x^{m+\frac12}_{h,\rho}|^2 \drho
\nonumber \\ & \qquad
+ \tfrac12 \left( e^{m+1}_{h,\rho}+ e^m_{h,\rho} , \eta_{h,\rho} \right) 
\nonumber \\ & \quad
= \left( \frac{d^{m+1}- d^m}{\Delta t} , \eta_h\, | x^{m+\frac12}_{h,\rho} |^2 
\right)
+ \left( x_t^{m+\frac12} - \frac{x^{m+1} - x^m}{\Delta t} , 
 \eta_h\, | x^{m+\frac12}_{h,\rho} |^2 \right) \nonumber \\ & \qquad
 + \left( x_t^{m+\frac12} , \eta_h\bigl( | x_\rho^{m+\frac12}|^2 -  | x^{m+\frac12}_{h,\rho} |^2  \bigr) \right) 
+ \left( x^{m+\frac12}_\rho - \tfrac12 ( x^{m+1}_\rho
+ x^m_\rho) , \eta_{h,\rho} \right)\nonumber \\ & \qquad
- \tfrac16 \sum_{j=1}^J h_j^2 \int_{I_j} \frac{x^{m+1}_{\rho}-x^m_{\rho}}{\Delta t} \cdot \eta_{h,\rho} \, |x^{m+\frac12}_{h,\rho}|^2 \drho
=: \sum_{\ell=1}^5 \langle S_\ell, \eta_h\rangle 
\qquad \forall\ \eta_h\in \Vh, \label{eq:err1}
\end{align} 
where we use the notation
$\langle S_\ell, \eta_h\rangle = S_\ell(\eta_h)$ for the specified linear
operators $S_\ell : \Vh \to \bR$, $\ell=1,\ldots,5$.
Choosing $\eta_h=e^{m+1}_h - e^m_h$ in \eqref{eq:err1} and noting
\eqref{eq:xhbound2} as well as \eqref{eq:decass}, we obtain
\begin{equation} \label{eq:rhsS}
\frac{c_0^2}{4 \Delta t} \| e^{m+1}_h - e^m_h \|^2_0 +\frac{c h^2}{\Delta t} \vert e^{m+1}_h - e^m_h \vert^2_1 +  \tfrac12 | e^{m+1}_{h} |^2_1 - \tfrac12 | e^m_{h} |^2_1 
\leq \sum_{\ell=1}^5 \langle S_\ell,e^{m+1}_h -e^m_h \rangle.
\end{equation}
Clearly, we have from \eqref{eq:xhbound2}, \eqref{eq:estIh} and the smoothness of $x$ that
\begin{align} \label{eq:boundS1}
| \langle S_1,e^{m+1}_h - e^m_h \rangle | & \leq c \sup_{t_m \leq t \leq t_{m+1}} \| x_t(\cdot,t) - I_h x_t(\cdot,t) \|_0 \| e^{m+1}_h - e^m_h \|_0 \nonumber\\
&\leq \frac{\epsilon}{\Delta t} \| e^{m+1}_h - e^m_h \|^2_0 + c_\epsilon h^4 \Delta t \sup_{t_m \leq t \leq t_{m+1}} \| x_t(\cdot,t) \|_2^2.
\end{align}
In addition, by Taylor expansion 
\begin{equation} \label{eq:taylortime}
 \left\| x_t^{m+\frac12} - \frac{x^{m+1} - x^m}{\Delta t} \right\|_{0,\infty} \leq c (\Delta t)^2,
 \end{equation}
 so that we obtain, again recalling \eqref{eq:xhbound2}, 
that 
\begin{equation} \label{eq:boundS2}
| \langle S_2,e^{m+1}_h - e^m_h \rangle | \leq c (\Delta t)^2 \| e^{m+1}_h - e^m_h \|_0 \leq  \frac{\epsilon}{\Delta t} \| e^{m+1}_h - e^m_h \|^2_0 + c_\epsilon (\Delta t)^5.
\end{equation}
Similarly, since $\|x^{m+\frac12}_{\rho \rho} - \tfrac12 ( x^{m+1}_{\rho \rho}+ x^m_{\rho \rho})\|_{0,\infty} \leq c(\Delta t)^2$, it holds that
\begin{align} \label{eq:boundS4}
| \langle S_4,e^{m+1}_h - e^m_h \rangle | & = 
\left| \left( x^{m+\frac12}_{\rho \rho} - \tfrac12 ( x^{m+1}_{\rho \rho}+ x^m_{\rho \rho})  , e^{m+1}_h - e^m_h  \right) \right| \nonumber \\ & 
\leq 
\frac{\epsilon}{\Delta t} \| e^{m+1}_h - e^m_h \|^2_0 + c_\epsilon (\Delta t)^5.
\end{align}
Next, let us write
\begin{align*}
\langle S_3,e^{m+1}_h - e^m_h \rangle & = 
\left( x_t^{m+\frac12} \cdot ( e^{m+1}_h - e^m_h) , | (I_h x^{m+\frac12})_\rho|^2 - | x^{m+\frac12}_{h,\rho} |^2 \right) \\ & \quad
+ \left( x_t^{m+\frac12} \cdot ( e^{m+1}_h - e^m_h) , | x^{m+\frac12}_{\rho} |^2 -  | (I_h x^{m+\frac12})_\rho|^2  \right) = S_{3,1}+S_{3,2}.
\end{align*}
Clearly, we have from  \eqref{eq:xhbound2} that
\[
| S_{3,1} | \leq c \| e^{m+1}_h -e^m_h \|_0 | e^{m+\frac12}_{h} |_1 \leq \frac{\epsilon}{\Delta t} \| e^{m+1}_h - e^m_h \|^2_0 + c_\epsilon \Delta t
 | e^{m+\frac12}_{h} |^2_1 .
\]
Let us remark that if we were to proceed in the same way for $S_{3,2}$ using \eqref{eq:estIh}, we would obtain the suboptimal bound $| S_{3,2} | \leq \frac{\epsilon}{\Delta t} \| e^{m+1}_h - e^m_h \|^2_0 + c_\epsilon  h^2 \Delta t$.
To avoid this, we write instead
\begin{align*}
S_{3,2} & = 
- \left( x_t^{m+\frac12} \cdot ( e^{m+1}_h - e^m_h) , | d^{m+\frac12}_\rho |^2 \right) 
+2 \left( x_t^{m+\frac12} \cdot ( e^{m+1}_h - e^m_h) , x^{m+\frac12}_\rho \cdot d^{m+\frac12}_\rho \right) \\ & 
= S^{(1)}_{3,2} + S^{(2)}_{3,2} 
\end{align*}
and derive with the help of \eqref{eq:estIh} 
\[
| S^{(1)}_{3,2} | \leq c \| d^{m+\frac12}_\rho \|_{0,\infty}  \| d^{m+\frac12}_\rho \|_0 \| e^{m+1}_h -e^m_h \|_0 \leq \frac{\epsilon}{\Delta t} \| e^{m+1}_h - e^m_h \|^2_0 + c_\epsilon h^4 \Delta t .
\]
Furthermore, integration by parts, on noting again
\eqref{eq:estIh}, yields that
\begin{align*}
S^{(2)}_{3,2} & = - 2 \left( x_{t,\rho}^{m+\frac12}  \cdot ( e^{m+1}_h - e^m_h) , 
x^{m+\frac12}_\rho \cdot d^{m+\frac12} \right) \\  & \quad
- 2 \left( x_t^{m+\frac12} \cdot ( e^{m+1}_h - e^m_h) , x^{m+\frac12}_{\rho \rho}
 \cdot d^{m+\frac12} \right) \\  & \quad
 - 2 \left( x_t^{m+\frac12} \cdot ( e^{m+1}_{h,\rho} - e^m_{h,\rho}) , x^{m+\frac12}_\rho \cdot d^{m+\frac12} \right) \\ & 
\leq c h^2 \| e^{m+1}_h -e^m_h \|_0 - 2 \left( x_t^{m+\frac12} \cdot ( e^{m+1}_{h,\rho} - e^m_{h,\rho}) , x^{m+\frac12}_\rho \cdot d^{m+\frac12} \right) \\ & 
= -2 \left(x_t^{m+\frac32} \cdot  e^{m+1}_{h,\rho} , x^{m+\frac32}_\rho \cdot d^{m+\frac32}\right) 
+ 2 \left( x_t^{m+\frac12} \cdot  e^{m}_{h,\rho} , x^{m+\frac12}_\rho \cdot d^{m+\frac12} \right) \\ & \quad
- 2 \left( e^{m+1}_{h,\rho} , (x^{m+\frac12}_\rho \cdot d^{m+\frac12}) x_t^{m+\frac12} - (x^{m+\frac32}_\rho \cdot d^{m+\frac32}) x_t^{m+\frac32} \right) + c h^2
 \| e^{m+1}_h -e^m_h \|_0 \\ & 
\leq - 2 \left(x_t^{m+\frac32} \cdot e^{m+1}_{h,\rho} , x^{m+\frac32}_\rho \cdot d^{m+\frac32}\right) 
+ 2 \left( x_t^{m+\frac12} \cdot e^{m}_{h,\rho} , x^{m+\frac12}_\rho \cdot d^{m+\frac12}\right) \\ & \quad 
+ c h^2 \Delta t | e^{m+1}_{h} |_1  + c h^2 \| e^{m+1}_h -e^m_h \|_0.
\end{align*}
In conclusion,
\begin{align} \label{eq:boundS3}
\langle S_3,e^{m+1}_h - e^m_h \rangle & 
\leq- 2 \left (x_t^{m+\frac32} \cdot  e^{m+1}_{h,\rho} , x^{m+\frac32}_\rho \cdot d^{m+\frac32}\right) 
+ 2 \left(x_t^{m+\frac12} \cdot e^{m}_{h,\rho} , x^{m+\frac12}_\rho \cdot d^{m+\frac12}\right) \nonumber \\ & \quad 
+ \frac{\epsilon}{\Delta t} \| e^{m+1}_h - e^m_h \|^2_0 + c_\epsilon \Delta t | e^{m+\frac12}_{h} |^2_1 + c_\epsilon h^4 \Delta t + c \Delta t | e^{m+1}_{h} |^2_1.
\end{align}
Finally, we have
\begin{align*}
& \langle S_5,e^{m+1}_h-e^m_h \rangle = \tfrac16 \sum_{j=1}^J h_j^2 \int_{I_j} \frac{x^{m+1}_\rho - x^m_\rho}{\Delta t}  \cdot (e^{m+1}_{h,\rho} - e^m_{h,\rho})  \bigl( | x^{m+\frac12}_\rho |^2 -  |x^{m+\frac12}_{h,\rho}|^2 \bigr) \drho
\\ & \qquad 
- \tfrac16 \sum_{j=1}^J h_j^2 \int_{I_j} \frac{x^{m+1}_{\rho}-x^m_{\rho}}{\Delta t} \cdot 
(e^{m+1}_{h,\rho} - e^m_{h,\rho}) \, |x^{m+\frac12}_{\rho}|^2 \drho =: S_{5,1} + S_{5,2}.
\end{align*}
Using the smoothness of $x$ and an interpolation estimate we have
\begin{align*}
| S_{5,1} | & \leq c h^2  | e^{m+1}_h - e^m_h |_1 \bigl( \vert e^{m+\frac12}_h \vert_1 + \vert d^{m+\frac12} \vert_1 \bigr) \\ & \leq \epsilon \frac{h^2}{\Delta t} | e^{m+1}_h - e^m_h |^2_1
+ c_\epsilon \Delta t \bigl( h^2 \vert e^{m+\frac12}_{h} \vert_1^2 + h^4 \bigr).
\end{align*}
Next, we have that
\begin{align*} 
S_{5,2} & = - \tfrac16 \sum_{j=1}^J h_j^2 \int_{I_j} \frac{x^{m+2}_{\rho}- x^{m+1}_\rho}{\Delta t} \cdot e^{m+1}_{h,\rho}   \,   |x^{m+\frac32}_{\rho}|^2 \drho
\\ & \quad
+ \tfrac16 \sum_{j=1}^J h_j^2 \int_{I_j} \frac{x^{m+1}_{\rho} - x^m_\rho}{\Delta t}  \cdot e^{m}_{h,\rho} \, |x^{m+\frac12}_{\rho}|^2 \drho \\ & \quad
+ \tfrac16 \sum_{j=1}^J h_j^2 \int_{I_j}  \bigl( |x^{m+\frac32}_{\rho}|^2 \, \frac{x^{m+2}_\rho - x^{m+1}_\rho}{\Delta t} - |x^{m+\frac12}_{\rho}|^2 \, \frac{x^{m+1}_\rho - x^m_\rho}{\Delta t} \bigr)  \cdot e^{m+1}_{h,\rho} \drho \\ &
\leq -\tfrac16 \sum_{j=1}^J h_j^2 \int_{I_j} \frac{x^{m+2}_{\rho}- x^{m+1}_\rho}{\Delta t} \cdot e^{m+1}_{h,\rho} \, |x^{m+\frac32}_{\rho}|^2 \drho \\ & \quad
+ \tfrac16 \sum_{j=1}^J h_j^2 \int_{I_j} \frac{x^{m+1}_{\rho} - x^m_\rho}{\Delta t}  \cdot e^{m}_{h,\rho} \, |x^{m+\frac12}_{\rho}|^2 \drho 
+ c h^2 \Delta t \vert e^{m+1}_{h} \vert_1.
\end{align*}
Combining the above bounds we deduce that
\begin{align} \label{eq:boundS5}
\langle S_5,e^{m+1}_h-e^m_h \rangle & \leq \epsilon \frac{h^2}{\Delta t} | e^{m+1}_h - e^m_h |^2_1 + c_\epsilon \Delta t \bigl( | e^{m+\frac12} |^2_1 + | e^{m+1} |^2_1 + h^4 \bigr) \nonumber  \\ & \quad
- \tfrac16 \sum_{j=1}^J h_j^2 \int_{I_j} \frac{x^{m+2}_{\rho}- x^{m+1}_\rho}{\Delta t} \cdot e^{m+1}_{h,\rho} \, |x^{m+\frac32}_{\rho}|^2 \drho 
\nonumber \\ & \quad
+ \tfrac16 \sum_{j=1}^J h_j^2 \int_{I_j} \frac{x^{m+1}_{\rho} - x^m_\rho}{\Delta t}  \cdot e^{m}_{h,\rho} \, |x^{m+\frac12}_{\rho}|^2 \drho. 
\end{align}
Choosing $\epsilon$ sufficiently small in \eqref{eq:boundS1}, 
\eqref{eq:boundS2}, \eqref{eq:boundS3}, \eqref{eq:boundS4}, and \eqref{eq:boundS5}
and recalling \eqref{eq:defEm}, we derive from \eqref{eq:rhsS} the bound \eqref{eq:lem2}. 
\end{proof}

\vspace{3mm}
We are now in position to complete the proof of Theorem~\ref{thm:main}. To do so, we use the bound \eqref{eq:predest} in \eqref{eq:lem2}, and recall the first estimate in \eqref{eq:Eequiv}, so that
\begin{align} \label{eq:dteh}
 \frac{c_0^2}{8 \Delta t} \| e^{m+1}_h - e^m_h \|^2_0 +  E^{m+1} 
& \leq E^m+  c \Delta t | e^{m+1}_{h} |^2_1 + c \Delta t | e^m_{h} |^2_1 + c \Delta t ( h^4 + (\Delta t)^4 ) \nonumber \\ &
\leq c \Delta t E^{m+1} + (1+c \Delta t) E^m  + c \Delta t ( h^4 + (\Delta t)^4 ).
\end{align}
If we choose $h_0$ so small that $\gamma h^\frac14_0 < \frac1{2c}$,
then $c \Delta t < \frac12$, and so
$(1 - c \Delta t)^{-1} \leq 1 + 2 c \Delta t$. Hence
we deduce from \eqref{eq:dteh} and the induction hypothesis \eqref{eq:induction} that
\begin{align*}
E^{m+1} & \leq (1 + 2 c \Delta t) \left [(1+ c \Delta t) E^m + c \Delta t ( h^4 + (\Delta t)^4 )\right] \\ &
\leq (1+ 4 c \Delta t) ( h^4 + (\Delta t)^4 ) e^{\mu t_m} + 2 c \Delta t ( h^4 + (\Delta t)^4 ) \\
& \leq
(1+ 6 c\Delta t) ( h^4 + (\Delta t)^4 ) e^{\mu t_m} 
\leq ( h^4 + (\Delta t)^4 ) e^{6c\Delta t} e^{\mu t_{m}}
= ( h^4 + (\Delta t)^4 ) e^{\mu t_{m+1}},
\end{align*}
if we choose $\mu = 6c$. Thus \eqref{eq:induction} holds for $m+1$ and hence for all $0 \leq m \leq M$. 
In conclusion, it follows from \eqref{eq:induction}, \eqref{eq:Eequiv} and
the first line in \eqref{eq:dteh} that
\begin{displaymath}
\max_{0 \leq m \leq M} | e^m_h |^2_1 + \sum_{m=0}^{M-1} \Delta t \, \left\| \frac{e^{m+1}_h - e^m_h}{\Delta t} \right\|^2_0  \leq c \bigl( h^4 + (\Delta t)^4 \bigr).
\end{displaymath}
Since $e^0_h=0$, we obtain
\begin{align*}
\Vert e^m_h \Vert_0^2 & = \sum_{k=0}^{m-1} \bigl( \Vert e^{k+1}_h \Vert_0^2 - \Vert e^k_h \Vert_0^2 \bigr) \leq \sum_{k=0}^{m-1} \Delta t  \left\| \frac{e^{k+1}_h - e^{k}_h}{\Delta t} \right\|_0 \bigl( \Vert e^{k+1}_h \Vert_0 + \Vert e^k_h \Vert_0 \bigr) \\
& \leq c \left( \sum_{k=1}^{M-1} \Delta t \left\| \frac{e^{k+1}_h - e^k_h}{\Delta t} \right\|^2_0 \right)^{\frac12} \max_{1 \leq k \leq M} \Vert e^k_h \Vert_0, \qquad 1 \leq m \leq M,
\end{align*}
which yields $\max_{1 \leq m \leq M} \Vert e^m_h \Vert_0 \leq c \bigl( h^2 + (\Delta t)^2 \bigr)$. The estimate \eqref{eq:l2h1} now follows with the help of
standard interpolation bounds. Finally, in order to derive \eqref{eq:timedif} we write
\begin{displaymath}
x_t^{m+\frac12} - \frac{x^{m+1}_h - x^m_h}{\Delta t} =  \left( x_t^{m+\frac12} - \frac{x^{m+1}- x^m}{\Delta t} \right) + \frac{d^{m+1}-d^m}{\Delta t} - \frac{e^{m+1}_h - e^m_h}{\Delta t}
\end{displaymath}
and use the above error bound together with an interpolation estimate and \eqref{eq:taylortime}. This concludes the proof of Theorem \ref{thm:main}.
\end{proof}

\setcounter{equation}{0}
\section{A second order scheme for curve diffusion} \label{sec:cd} 

Similarly to the approach in Section~\ref{sec:csf}, here we aim to apply 
the predictor--corrector idea from \cite{JiangSZZ25} to the
first order method \eqref{eq:cdfea} from \cite{cd}. Due to the presence of the
curvature variable $y^m_h$ on the right hand side of \eqref{eq:cdfeaa},
special care must be taken in the preparation of the predictor and
corrector steps.
Simply solving \eqref{eq:cdfeaa} with $\Delta t$ replaced by
$\frac12 \Delta t$ for the predictor, and using the obtained solutions
for the explicit terms in \eqref{eq:cdfea} for the corrector step 
may lead to oscillations in the polygonal approximations.
However, the following scheme works well in practice.

For $m \geq 0$, given $(x^m_h,y^{m-\frac12}_h) \in \Vh\times\Vh$ first find 
$(x^{m+\frac12}_h, z^{m+\frac12}_h) \in \Vh \times \Vh$ such that
\begin{subequations} \label{eq:cdpredictor12}
\begin{align}
& \left( \frac{x^{m+\frac12}_h-x^m_h}{\frac12\Delta t} , \chi_h |x^{m}_{h,\rho}|^2 \right) 
- \left( z^{m+\frac12}_{h,\rho} , \chi_{h,\rho} \right) 
= 2 \left( (z^{m+\frac12}_{h,\rho} \cdot x^{m}_{h,\rho}) y^{m-\frac12}_h , \chi_h \right)
\nonumber \\ & \quad
+ \left( (y^{m-\frac12}_h \cdot z^{m+\frac12}_h) y^{m-\frac12}_h , \chi_h | x^{m}_{h,\rho} |^2 \right)
+  \left( F_2(x^{m}_{h,\rho},y^{m-\frac12}_h,y^{m-\frac12}_{h,\rho})z^{m+\frac12}_h , \chi_h \right) 
, \label{eq:cdpredictor12a} \\
& \left( z^{m+\frac12}_h , \eta_h |x^{m}_{h,\rho}|^2 \right) 
+ 2 \left(((x^{m+\frac12}_{h,\rho}-x^m_{h,\rho}) \cdot x^m_{h,\rho})
y^{m-\frac12}_h ,\eta_h \right)
+ \left( x^{m+\frac12}_{h,\rho} , \eta_{h,\rho} \right) = 0
\label{eq:cdpredictor12b}
\end{align}
\end{subequations}
for all $\chi_h, \eta_h \in \Vh$.
Then find $(x^{m+1}_h, y^{m+\frac12}_h) \in \Vh \times \Vh$ such that
\begin{subequations} \label{eq:cnprefd}
\begin{align}
& \left( \frac{x^{m+1}_h-x^m_h}{\Delta t} , \chi_h |x^{m+\frac12}_{h,\rho}|^2 \right) 
- \left( y^{m+\frac12}_{h,\rho} , \chi_{h,\rho} \right) 
= 2 \left( (y^{m+\frac12}_{h,\rho} \cdot x^{m+\frac12}_{h,\rho}) z^{m+\frac12}_h , \chi_h \right)
\nonumber \\ & \quad
+ \left( (z^{m+\frac12}_h \cdot y^{m+\frac12}_h) z^{m+\frac12}_h , \chi_h | x^{m+\frac12}_{h,\rho} |^2 \right)
+  \left( F_2(x^{m+\frac12}_{h,\rho},z^{m+\frac12}_h,z^{m+\frac12}_{h,\rho})y^{m+\frac12}_h , \chi_h \right) 
, \label{eq:cnprefda} \\
& \left( y^{m+\frac12}_h , \eta_h |x^{m+\frac12}_{h,\rho}|^2 \right) 
+ \tfrac12\left( x^{m+1}_{h,\rho} , \eta_{h,\rho} \right) =
- \tfrac12\left( x^{m}_{h,\rho} , \eta_{h,\rho} \right) 
\label{eq:cnprefdb}
\end{align}
\end{subequations}
for all $\chi_h, \eta_h \in \Vh$.
Observe that the second term in \eqref{eq:cdpredictor12b} is motivated by the 
linearization
$| x^{m+\frac12}_{h,\rho} |^2 \approx | x^m_{h,\rho} |^2 +
2 (x^m_{h,\rho},x^{m+\frac12}_{h,\rho}-x^m_{h,\rho})$.
Moreover, $z^{m+\frac12}_h$ acts as a predictor for $y^{m+\frac12}_h$, and both
approximate the true solution $y^{m+\frac12}$.
For the purposes of the predictor
step it is sufficiently accurate to make use of $y^{m-\frac12}_h$,
even though it is shifted by half a time step from $t_m$.

Regardless of how the predictor variables $(x^{m+\frac12}_h, z^{m+\frac12}_h) \in \Vh \times \Vh$ are obtained, the following theorem proves existence, uniqueness and unconditional stability of the corrector step \eqref{eq:cnprefd}. 

\begin{theorem} 
Assume that $|x^{m+\frac12}_{h,\rho}| > 0$ in I. Then
there exists a unique solution
$(x^{m+1}_h, y^{m+\frac12}_h) \in \Vh \times \Vh$ to \eqref{eq:cnprefd},
satisfying the stability estimate
\begin{equation} \label{eq:cdstab}
\tfrac12 |x^{m+1}_h|^2_1 + 
\Delta t \left( | y^{m+\frac12}_{h,\rho} + 
(z^{m+\frac12}_h \cdot y^{m+\frac12}_h) x^{m+\frac12}_{h,\rho} |^2, 1 \right) 
= \tfrac12 |x^m_h|^2_1.
\end{equation}
\end{theorem}
\begin{proof}
The well-posedness of the linear system \eqref{eq:cdfea} has been shown in 
\cite[Theorem~5.1]{cd}. The proof for \eqref{eq:cnprefd} is nearly identical,
and so we omit it here.

In order to prove \eqref{eq:cdstab}, we choose 
$\chi = y^{m+\frac12}_h$ in \eqref{eq:cnprefda} and $\eta = x^{m+1}_h - x^m_h$ in
\eqref{eq:cnprefdb} to yield that 
\begin{align*}
\tfrac12 |x^{m+1}_h|^2_1 - \tfrac12 |x^m_h|^2_1
& = \tfrac12 \left( x^{m+1}_{h,\rho} + x^m_{h,\rho} , x^{m+1}_{h,\rho} - x^m_{h,\rho} \right)
\nonumber \\ &
= - \Delta t |y^{m+\frac12}_h|^2_1
- 2 \Delta t \left( y^{m+\frac12}_{h,\rho} \cdot x^{m+\frac12}_{h,\rho} ,
z^{m+\frac12}_h \cdot y^{m+\frac12}_h \right) \nonumber \\ & \quad
- \Delta t \left( (z^{m+\frac12}_h \cdot y^{m+\frac12}_h)^2 , | x^{m+\frac12}_{h,\rho} |^2 \right)
\nonumber \\ &
= - \Delta t \left( | y^{m+\frac12}_{h,\rho} + (z^{m+\frac12}_h \cdot y^{m+\frac12}_h) 
x^{m+\frac12}_{h,\rho} |^2, 1 \right),
\end{align*}
where we have used that $F_2(a,b,c)$ is anti-symmetric, recall 
\eqref{eq:defF2}.
\end{proof}

\begin{remark} 
An alternative predictor step, which also works well in practice, is given
by the following fully implicit scheme.
Find $(x^{m+\frac12}_h, z^{m+\frac12}_h) \in \Vh \times \Vh$ such that
\begin{subequations} \label{eq:cdpredictornl}
\begin{align}
& \left( \frac{x^{m+\frac12}_h-x^m_h}{\frac12\Delta t} , \chi_h |x^{m+\frac12}_{h,\rho}|^2 \right) 
- \left( z^{m+\frac12}_{h,\rho} , \chi_{h,\rho} \right) 
= 2 \left( (z^{m+\frac12}_{h,\rho} \cdot x^{m+\frac12}_{h,\rho}) z^{m+\frac12}_h , \chi_h \right)
\nonumber \\ & \quad
+ \left( (z^{m+\frac12}_h \cdot z^{m+\frac12}_h) z^{m+\frac12}_h , \chi_h | x^{m+\frac12}_{h,\rho} |^2 \right)
+  \left( F_2(x^{m+\frac12}_{h,\rho},z^{m+\frac12}_h,z^{m+\frac12}_{h,\rho})z^{m+\frac12}_h , \chi_h \right) 
, \label{eq:cdpredictornla} \\
& \left( z^{m+\frac12}_h , \eta_h |x^{m+\frac12}_{h,\rho}|^2 \right) 
+ \left( x^{m+\frac12}_{h,\rho} , \eta_{h,\rho} \right) = 0
\label{eq:cdpredictornlb}
\end{align}
\end{subequations}
for all $\chi_h, \eta_h \in \Vh$.
Compared to \eqref{eq:cdpredictor12}, the above is slightly more natural, more
compact and does not require knowledge of $y^{m-\frac12}_h$. It also provides
the motivation for the linearized step \eqref{eq:cdpredictor12}.
But as \eqref{eq:cdpredictornl} requires the solution of a nonlinear system of
equations at each time step, we prefer the linear predictor step
\eqref{eq:cdpredictor12} for the numerical simulations in this paper.
\end{remark}

\setcounter{equation}{0}
\section{Numerical results} \label{sec:nr}

We implemented all our schemes within the
finite element toolbox Alberta, \cite{Alberta}, using
the sparse factorization package UMFPACK, see \cite{Davis04},
for the solution of the linear systems of equations arising at each time 
level.

For all our numerical simulations we use a uniform partitioning of $[0,1]$,
so that $q_j = jh$, $j=0,\ldots,J$, with $h = \frac 1J$. 
For the initial data we choose $x_h^0 = I_h x_0$.
In addition, for the schemes for curve diffusion, we let $y^0_h \in \Vh$ 
be the solution of 
\begin{equation*}
\left( y^{0}_h , \eta_h |x^0_{h,\rho}|^2 \right) 
+ \left( x^{0}_{h,\rho} , \eta_{h,\rho} \right) = 0
\qquad \forall\ \eta_h \in \Vh,
\end{equation*}
compare with \eqref{eq:cdfeab}, and set $y^{-\frac12}_h = y^0_h$.
For some of our simulations we will monitor the ratio
\begin{equation} \label{eq:ratio}
\ratio^m = \dfrac{\max_{j=1,\ldots,J} |x_h^m(q_j) - x_h^m(q_{j-1})|}
{\min_{j=1,\ldots, J} |x_h^m(q_j) - x_h^m(q_{j-1})|}
\end{equation}
between the lengths of the longest and shortest element of $\Gamma^m=x^m_h(I)$.
Clearly $\ratio^m\geq1$, with equality if and only if the curve is 
equidistributed.

\newcommand{\errorxL}{\| x -  x_h\|_0}
\newcommand{\errorxH}{\| x -  x_h\|_1}
\newcommand{\erroryL}{\| y -  y_h\|_0}
\newcommand{\erroryH}{\| y -  y_h\|_1}

\subsection{Curve shortening flow}

For a convergence experiment we use the well known solution of a family of 
shrinking circles with radius
$r(t) = (r^2(0) - 2t)^\frac12$.
To this end, we choose the particular parameterization
\begin{equation} \label{eq:true}
x(\rho,t)=\hat x(\rho,t)=(1-2t)^{\frac12} \binom{\cos g(\rho)}{\sin g(\rho)},
\quad \text{where }\ 
g(\rho) = 2\pi\rho + \delta \sin(2\pi\rho), \ \delta = 0.1,
\end{equation}
and define
\begin{equation*} 
f_{csf} = |\hat x_\rho|^2 \hat x_t - \hat x_{\rho\rho}
\end{equation*}
as an extra forcing for \eqref{eq:DD95}.
In particular, we use it to add
$(f_{csf}(\cdot, t_{m+\frac12}) , \eta_h )^h$
to the right hand sides of \eqref{eq:pred} and \eqref{eq:csfd}, respectively. 
We also define the errors 
\begin{equation} \label{eq:errorxL}
\errorxL = \max_{m=0,\ldots,M} \| x(\cdot,t_m) -  x^m_h\|_0,\quad
\errorxH = \max_{m=0,\ldots,M} \| x(\cdot,t_m) -  x^m_h\|_1.
\end{equation}
See Table~\ref{tab:csfcn2_delta01}, where the optimal convergence rates
from Theorem~\ref{thm:main} are confirmed.
Here we partition the time interval $[0,T]$, with $T=0.25$, into uniform
time steps of size $\Delta t = h$, for $h = J^{-1} = 2^{-k}$, 
$k=5,\ldots,12$. 
In the same table we also report on the errors for the scheme
\eqref{eq:DD95:fea} from \cite{DeckelnickD95}, 
which demonstrates only first order convergence.
The fact that the $H^1$--errors for \eqref{eq:DD95:fea} and \eqref{eq:fea} are
identical to the displayed digits suggests that the errors from the 
spatial discretization are dominating for that norm.
\begin{table}
\center
\begin{tabular}{|r|c|c|c|c|c|c|c|c|}
\hline
& \multicolumn{4}{c|}{\eqref{eq:DD95:fea}} & \multicolumn{4}{c|}{\eqref{eq:fea}} \\
$J$ & $\errorxL$ & EOC & $\errorxH$ & EOC & $\errorxL$ & EOC & $\errorxH$ & EOC \\ \hline
32   &1.7810e-02& ---&3.6211e-01& ---&3.2653e-03& ---&3.6211e-01& --- \\
64   &1.0145e-02&0.81&1.8113e-01&1.00&8.1694e-04&2.00&1.8113e-01&1.00 \\
128  &5.4018e-03&0.91&9.0578e-02&1.00&2.0427e-04&2.00&9.0578e-02&1.00 \\
256  &2.7861e-03&0.96&4.5290e-02&1.00&5.1070e-05&2.00&4.5290e-02&1.00 \\
512  &1.4147e-03&0.98&2.2645e-02&1.00&1.2768e-05&2.00&2.2645e-02&1.00 \\
1024 &7.1285e-04&0.99&1.1323e-02&1.00&3.1920e-06&2.00&1.1323e-02&1.00 \\
2048 &3.5780e-04&0.99&5.6613e-03&1.00&7.9799e-07&2.00&5.6613e-03&1.00 \\
4096 &1.7925e-04&1.00&2.8307e-03&1.00&1.9950e-07&2.00&2.8307e-03&1.00 \\
\hline
\end{tabular}
\caption{Errors for the convergence test for \eqref{eq:true} 
over the time interval $[0,0.25]$.
We also display the experimental orders of convergence (EOC).}
\label{tab:csfcn2_delta01}
\end{table}%
\begin{table}
\center
\begin{tabular}{|r|c|r|c|r|}
\hline
& \multicolumn{2}{c|}{\eqref{eq:DD95:fea}} & \multicolumn{2}{c|}{\eqref{eq:fea}} \\
$-\log_2 \Delta t$ & $\errorxL$ & CPU time [s] & $\errorxL$ & CPU time [s] \\ 
\hline
 9 & 1.4235e-03 & 2   &2.5407e-06 & 4 \\
10 & 7.1493e-04 & 4   &5.3356e-07 & 8 \\
11 & 3.5822e-04 & 9   &1.9950e-07 & 17\\
12 & 1.7925e-04 & 17  &1.9950e-07 & 33\\
13 & 8.9605e-05 & 34  &1.9950e-07 & 67\\
14 & 4.4746e-05 & 68  &1.9950e-07 & 133\\
15 & 2.2307e-05 & 137 &1.9950e-07 & 269\\
16 & 1.1085e-05 & 274 &1.9950e-07 & 535\\
17 & 5.4732e-06 & 548 &1.9950e-07 &1072\\
18 & 2.6673e-06 &1102 &1.9950e-07 &2159\\
19 & 1.2644e-06 &2193 &1.9950e-07 &4283\\ 
20 & 5.6319e-07 &4403 &1.9950e-07 &8567\\
\hline
\end{tabular}
\caption{Errors and CPU times for the convergence test for \eqref{eq:true} 
over the time interval $[0,0.25]$. Here $J=4096$ is fixed.
}
\label{tab:csfCPU}
\end{table}%
The dramatic computational superiority of the second order scheme becomes
apparent in Table~\ref{tab:csfCPU}. Here we compare our scheme \eqref{eq:fea} 
once again with the first order scheme \eqref{eq:DD95:fea}. This time the
spatial discretization is fixed, with $J=4096$, while we choose 
$\Delta t = 2^{-k}$, $k=9,\ldots,20$. We note that already for $k=11$ the spatial
error dominates for the scheme \eqref{eq:fea}, recall \eqref{eq:l2h1}. 
Moreover, to reach the same error as the second order scheme \eqref{eq:fea} 
with $\Delta t = 2^{-10}$, the first order scheme \eqref{eq:DD95:fea} needs
to employ a time step size of $\Delta t = 2^{-20}$, which results in a CPU time
that is 550 times larger.

\begin{figure}
\center
\includegraphics[angle=-90,width=0.3\textwidth]{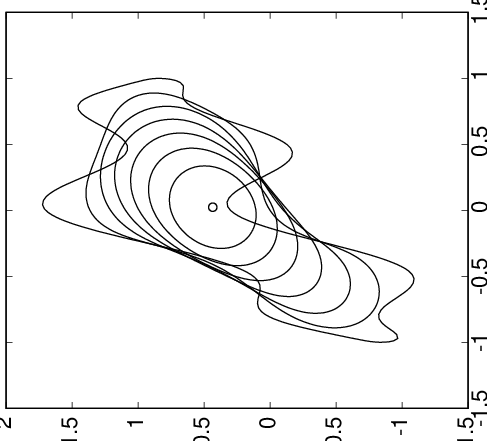}\\
\includegraphics[angle=-90,width=0.31\textwidth]{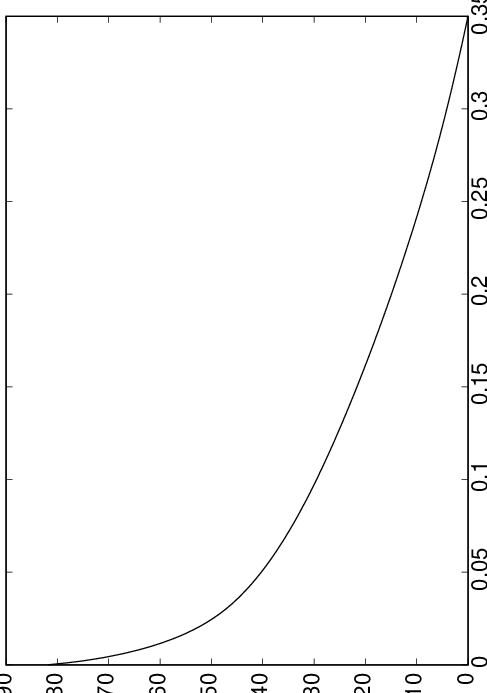}
\includegraphics[angle=-90,width=0.31\textwidth]{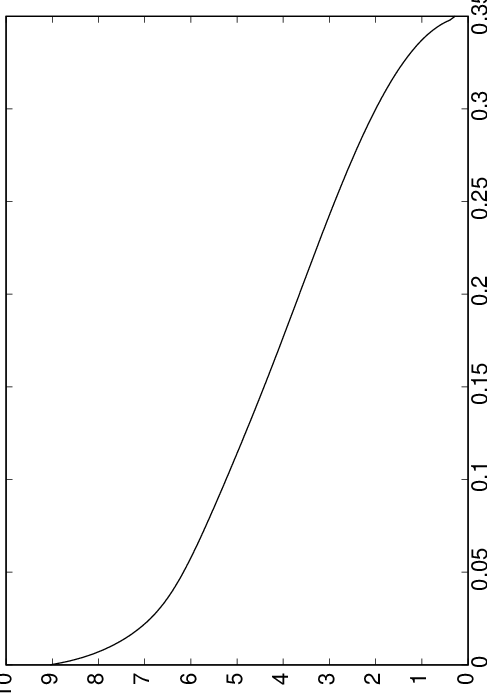}
\includegraphics[angle=-90,width=0.31\textwidth]{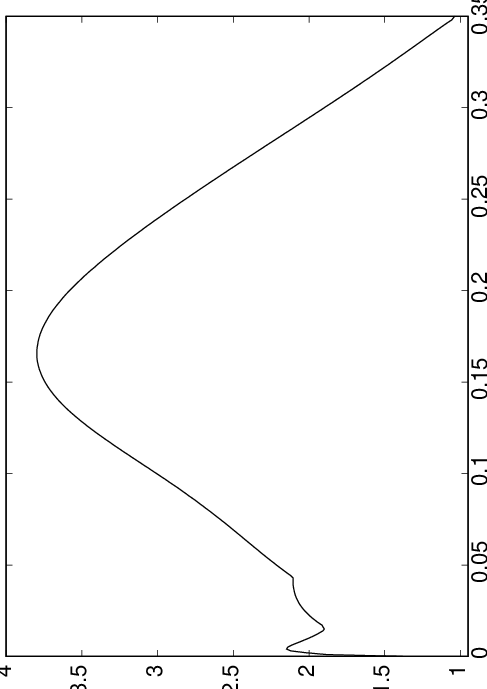}\\
{\footnotesize
\hspace{2cm} $\int_I |x^m_{h,\rho}|^2\drho$ 
\hspace{3.5cm} $|\Gamma^m|$ 
\hspace{4.2cm} $\ratio^m$} \hspace{2cm}
\caption{Curve shortening flow for the nonconvex curve \eqref{eq:mikula}. 
On top we show $\Gamma^m$ at times $t=0,0.05,\ldots,T=0.35$.
Below we show the evolutions of $\int_I |x^m_{h,\rho}|^2\drho$ (left),
$|\Gamma^m|$ (middle) and $\ratio^m$ (right) over time.
}
\label{fig:mikula}
\end{figure}

We also present the evolution for an initially nonconvex planar curve under
curve shortening flow. To this end, we let
\begin{equation} \label{eq:mikula}
 x(\rho,0)=\binom{\cos{u(\rho)}}
  {\tfrac12\sin{u(\rho)} + \sin{(\cos{u(\rho)})} + \sin{u(\rho)}\,
  [\tfrac15+\sin{u(\rho)}\,\sin^2{u(3\rho)}]}, 
\end{equation}
where $u(\rho) = 2\pi\rho$, see also \cite[Fig.\ 1]{MikulaS01}.
The evolution of the curves $\Gamma^m=x^m_h(I)$ for the solutions of the scheme
\eqref{eq:fea} is shown in Figure~\ref{fig:mikula}. Here we employ the
discretization parameters $J=256$ and $\Delta t = 10^{-3}$.
We observe that the curve first becomes
convex, and then shrinks to a point. Both the Dirichlet energy of $x^m$ and the
length of the curve $\Gamma^m$ are monotonically decreasing, while the ratio
\eqref{eq:ratio} at first increases slightly to about 3.8, before approaching
unity towards the end of the simulation.

\subsection{Curve diffusion}

We begin with a convergence experiment for the prescribed solution of a 
translated and dilated circle parameterized by
\begin{equation} \label{eq:solx}
x(\rho,t)=\hat x(\rho,t)=\binom{t^2}{t^2} + (1 + t^3) \binom{\cos g(\rho)}{\sin g(\rho)},
\end{equation}
recall the definition of $g$ from \eqref{eq:true},
by constructing the right-hand side 
\[
f_{cd} = |\hat x_\rho|^2 \hat x_t + \hat y_{\rho\rho} - F_{cd}(\hat x_\rho,\hat y,\hat y_\rho) \hat y, \quad \text{where } \hat y = \frac{\hat x_{\rho\rho}}{|\hat x_\rho|^2}, 
\]
for \eqref{eq:xt2}.
Upon adding the correction term
$( f_{cd}(\cdot, t_{m+\frac12}) , \eta_h )^h$
to the right hand sides of \eqref{eq:cdpredictor12a} 
and \eqref{eq:cnprefda}, we can perform a convergence
experiment for the scheme \eqref{eq:cdpredictor12},
\eqref{eq:cnprefd}, comparing the obtained discrete solutions with
\eqref{eq:solx}. In particular, we define the errors $\erroryL$ and $\erroryH$
similarly to \eqref{eq:errorxL}, upon defining
$y^{m+1}_h = 2 y^{m+\frac12}_h - y^m_h$ for $m=0,\ldots,M-1$.
The obtained errors are displayed in Table~\ref{tab:cn7cxy},
where we can see second order convergence. We compare the errors with the ones
from the scheme \eqref{eq:cdfea} in Table~\ref{tab:cn0xy}, showing only
first order convergence.
In both cases we partition the time interval $[0,T]$, with $T=1$, into uniform
time steps of size $\Delta t = h$, for $h = J^{-1} = 2^{-k}$, 
$k=5,\ldots,12$. 
\begin{table}
\center
\begin{tabular}{|r|c|c|c|c|c|c|c|c|}
\hline
$J$ & $\errorxL$ & EOC  & $\errorxH$ & EOC  & $\erroryL$ & EOC  & $\erroryH$ & EOC \\ \hline
32   &3.8513e-03& ---&7.2497e-01& ---&3.3870e-03& ---&3.8417e-01& --- \\
64   &9.4550e-04&2.03&3.6237e-01&1.00&8.4367e-04&2.01&1.9102e-01&1.01 \\
128  &2.3453e-04&2.01&1.8117e-01&1.00&2.1074e-04&2.00&9.5376e-02&1.00 \\
256  &5.8423e-05&2.01&9.0582e-02&1.00&5.2673e-05&2.00&4.7671e-02&1.00 \\
512  &1.4581e-05&2.00&4.5291e-02&1.00&1.3167e-05&2.00&2.3834e-02&1.00 \\
1024 &3.6421e-06&2.00&2.2645e-02&1.00&3.2918e-06&2.00&1.1917e-02&1.00 \\
2048 &9.1015e-07&2.00&1.1323e-02&1.00&8.2295e-07&2.00&5.9582e-03&1.00 \\
4096 &2.2749e-07&2.00&5.6613e-03&1.00&2.0574e-07&2.00&2.9791e-03&1.00 \\
\hline
\end{tabular}
\caption{(\eqref{eq:cdpredictor12}, \eqref{eq:cnprefd})
Errors for the convergence test for \eqref{eq:solx} 
over the time interval $[0,1]$.
We also display the experimental orders of convergence (EOC).}
\label{tab:cn7cxy}
\end{table}%
\begin{table}
\center
\begin{tabular}{|r|c|c|c|c|c|c|c|c|}
\hline
$J$ & $\errorxL$ & EOC  & $\errorxH$ & EOC  & $\erroryL$ & EOC  & $\erroryH$ & EOC \\ \hline
32   &2.7927e-02& ---&7.3353e-01& ---&6.3618e-02& ---&4.9977e-01& --- \\
64   &1.3609e-02&1.04&3.6833e-01&0.99&3.1180e-02&1.03&2.4354e-01&1.04 \\
128  &6.7142e-03&1.02&1.8458e-01&1.00&1.5441e-02&1.01&1.2029e-01&1.02 \\
256  &3.3342e-03&1.01&9.2399e-02&1.00&7.6840e-03&1.01&5.9790e-02&1.01 \\
512  &1.6613e-03&1.01&4.6227e-02&1.00&3.8330e-03&1.00&2.9807e-02&1.00 \\
1024 &8.2922e-04&1.00&2.3120e-02&1.00&1.9143e-03&1.00&1.4882e-02&1.00 \\
2048 &4.1425e-04&1.00&1.1562e-02&1.00&9.5658e-04&1.00&7.4354e-03&1.00 \\
4096 &2.0703e-04&1.00&5.7814e-03&1.00&4.7815e-04&1.00&3.7163e-03&1.00 \\
\hline
\end{tabular}
\caption{(\eqref{eq:cdfea})
Errors for the convergence test for \eqref{eq:solx} 
over the time interval $[0,1]$.
We also display the experimental orders of convergence (EOC).}
\label{tab:cn0xy}
\end{table}%
\begin{table}
\center
\begin{tabular}{|r|c|r|c|r|}
\hline
& \multicolumn{2}{c|}{\eqref{eq:cdfea}} & \multicolumn{2}{c|}{\eqref{eq:cdpredictor12}, \eqref{eq:cnprefd}} \\
$-\log_2 \Delta t$ & $\errorxL$ & CPU time [s] & $\errorxL$ & CPU time [s] \\ 
\hline
 7 & 6.6487e-03 & 3  & 1.3428e-04 & 7   \\
 8 & 3.3177e-03 & 6  & 3.3432e-05 & 14  \\
 9 & 1.6572e-03 & 12 & 8.3974e-06 & 28  \\
10 & 8.2824e-04 & 25 & 2.1616e-06 & 56  \\
11 & 4.1405e-04 & 50 & 6.0845e-07 & 113 \\
12 & 2.0703e-04 & 99 & 2.2749e-07 & 222 \\
13 & 1.0354e-04 & 198& 2.0075e-07 & 449 \\
14 & 5.1803e-05 & 396& 2.0049e-07 & 886 \\
15 & 2.5934e-05 & 792& 2.0043e-07 & 1813\\
16 & 1.3000e-05 &1592& 2.0041e-07 & 3503\\ 
17 & 6.5330e-06 &3156& 2.0041e-07 & 7135\\
18 & 3.3001e-06 &6326& 2.0041e-07 & 13932\\
19 & 1.6847e-06 &12551&2.0041e-07 & 27294\\
20 & 8.7857e-07 &25108&2.0041e-07 & 53700\\
\hline
\end{tabular}
\caption{Errors and CPU times for the convergence test for \eqref{eq:solx} 
over the time interval $[0,1]$. Here $J=4096$ is fixed.
}
\label{tab:cdCPU}
\end{table}%
Similarly to Table~\ref{tab:csfCPU}, we once again see the superiority of the 
second order scheme when we keep $J$ fixed in Table~\ref{tab:cdCPU}. In
particular, to achieve an overall $L^2$--error of about $2\times10^{-6}$, 
the first order method \eqref{eq:cdfea} needs about 220 times longer than the 
second order method \eqref{eq:cdpredictor12}, \eqref{eq:cnprefd}.

\begin{figure}
\center
\includegraphics[angle=-0,width=0.3\textwidth]{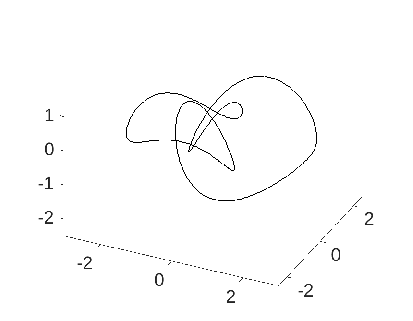}
\includegraphics[angle=-0,width=0.3\textwidth]{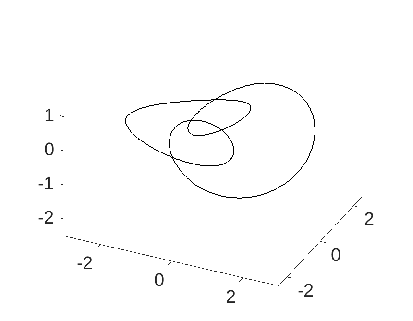}
\includegraphics[angle=-0,width=0.3\textwidth]{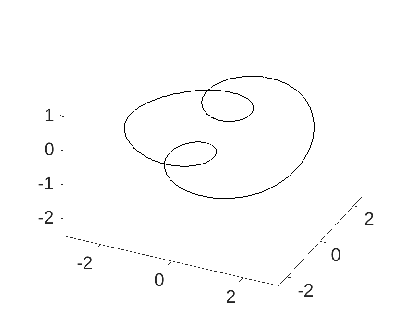}
\includegraphics[angle=-0,width=0.3\textwidth]{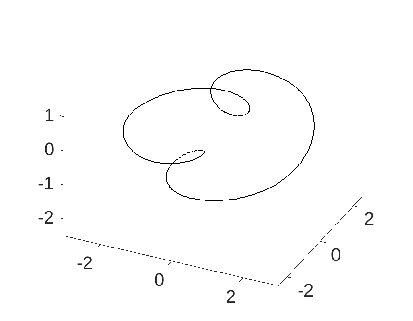}
\includegraphics[angle=-0,width=0.3\textwidth]{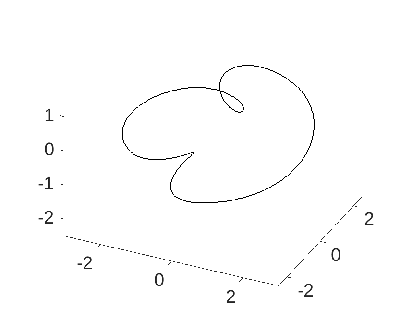}
\includegraphics[angle=-0,width=0.3\textwidth]{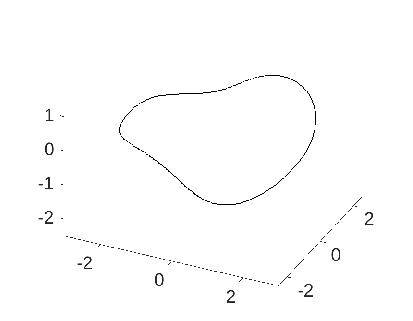}
\includegraphics[angle=-0,width=0.3\textwidth]{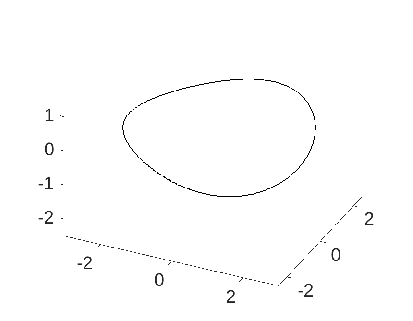}
\includegraphics[angle=-0,width=0.3\textwidth]{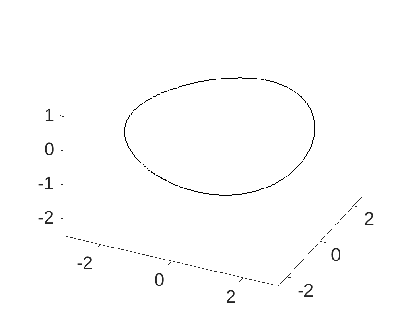}
\includegraphics[angle=-0,width=0.3\textwidth]{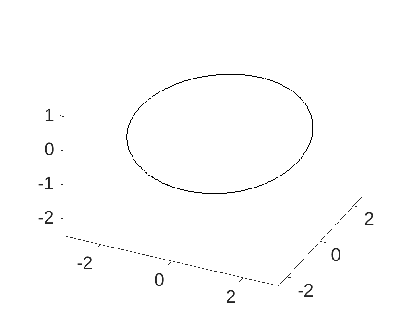} \\
\includegraphics[angle=-90,width=0.3\textwidth]{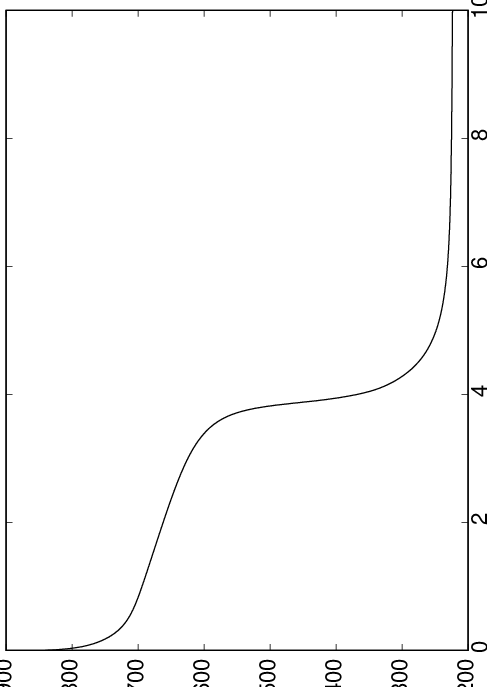}
\includegraphics[angle=-90,width=0.3\textwidth]{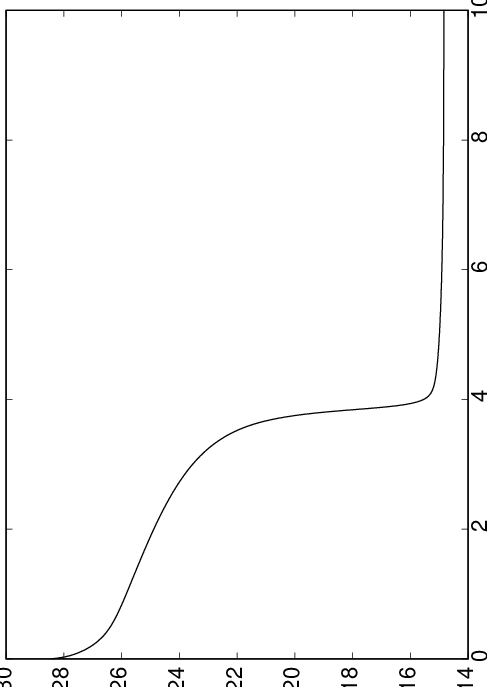}
\includegraphics[angle=-90,width=0.3\textwidth]{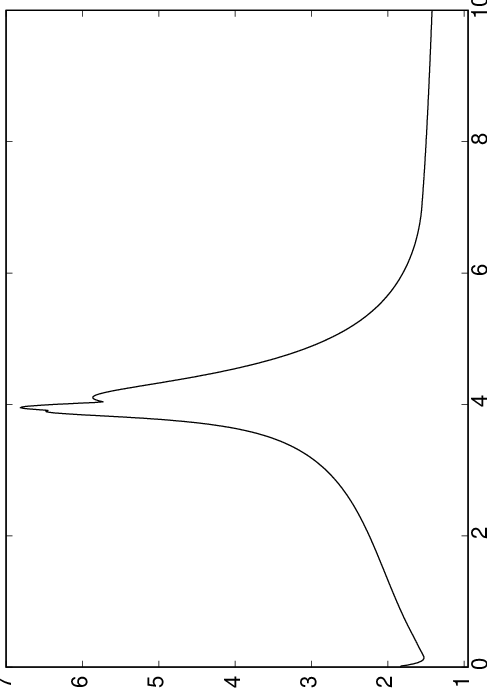}\\
{\footnotesize
\hspace{2cm} $\int_I |x^m_{h,\rho}|^2\drho$ 
\hspace{3.5cm} $|\Gamma^m|$ 
\hspace{4.2cm} $\ratio^m$} \hspace{2cm}
\caption{Curve diffusion for the two interlocked rings \eqref{eq:irings}. 
We show $\Gamma^m$ at times $t= 0, 1, 3, 3.5, 3.7, 4, 5, 6, T=10$.
Below we show the evolutions of $\int_I |x^m_{h,\rho}|^2\drho$ (left),
$|\Gamma^m|$ (middle) and $\ratio^m$ (right) over time.
}
\label{fig:irings}
\end{figure}

We also present a numerical experiment for $d=3$. To this end, we consider
the two interlocked rings in $\bR^3$ from \cite[Fig.\ 4]{cd}. 
In particular, the initial curve is given by
\begin{equation}
x_0(\rho) = \tfrac18 
\begin{pmatrix}
10 (\cos(2\pi\rho)+\cos(6\pi\rho))+\cos(4\pi\rho)
+\cos(8\pi\rho) \\
6\sin(2\pi\rho)+10\sin(6\pi\rho) \\
4\sin(6\pi\rho)\sin(5\pi\rho)+4\sin(8\pi\rho)-2\sin(12\pi\rho)
\end{pmatrix}
\quad \rho \in I.
\label{eq:irings}
\end{equation}
The evolution of this curve under curve diffusion can
be seen in Figure~\ref{fig:irings}, showing perfect agreement with 
\cite[Fig.\ 4]{cd}. But here, for the scheme \eqref{eq:cdpredictor12},
\eqref{eq:cnprefd}, we employ the
discretization parameters $J=512$ and $\Delta t = 10^{-2}$, i.e.\ a time step
size that is 100 times larger than in \cite{cd}.
As a consequence, the computation took only 9 seconds, compared to the
421 seconds for the simulation for the first order scheme \eqref{eq:cdfea} 
in \cite[Fig.\ 4]{cd}.
We note once again that the ratio \eqref{eq:ratio} remains bounded during this
simulation, approaching a value close to unity as the evolution settles on the
circular steady state solution.

\def\soft#1{\leavevmode\setbox0=\hbox{h}\dimen7=\ht0\advance \dimen7
  by-1ex\relax\if t#1\relax\rlap{\raise.6\dimen7
  \hbox{\kern.3ex\char'47}}#1\relax\else\if T#1\relax
  \rlap{\raise.5\dimen7\hbox{\kern1.3ex\char'47}}#1\relax \else\if
  d#1\relax\rlap{\raise.5\dimen7\hbox{\kern.9ex \char'47}}#1\relax\else\if
  D#1\relax\rlap{\raise.5\dimen7 \hbox{\kern1.4ex\char'47}}#1\relax\else\if
  l#1\relax \rlap{\raise.5\dimen7\hbox{\kern.4ex\char'47}}#1\relax \else\if
  L#1\relax\rlap{\raise.5\dimen7\hbox{\kern.7ex
  \char'47}}#1\relax\else\message{accent \string\soft \space #1 not
  defined!}#1\relax\fi\fi\fi\fi\fi\fi}

\end{document}